%{
%{
\font\de=cmssi12
\font\twelverm=cmr12

% This file defines the standard font setup, including the macros
% \tenpoint and \twelvepoint, as well as some standard dimension settings.
%
% It is indended to be used as a basis for most other ``personal'' formats.

% Make sure we haven't already been run

\ifx\UsualIsLoaded\undefined
\let\UsualIsLoaded=\relax		% define UsualIsLoaded

\font\fourteenrm=cmr12  scaled \magstep1
\font\fourteenbf=cmbx12 scaled \magstep1
\font\fourteentt=cmtt12 scaled \magstep1
\font\fourteensl=cmsl12 scaled \magstep1
\font\fourteensy=cmsy10 scaled \magstep2
\font\fourteeni=cmmi12  scaled \magstep1
\font\fourteenit=cmti12 scaled \magstep1
\font\fourteensc=cmcsc10 scaled \magstep2
\font\fourteenbit=cmssi12 scaled \magstep1
\font\fourteenbbb=msym10 scaled \magstep2

\font\twelverm=cmr12
\font\twelvebf=cmbx12
\font\twelvett=cmtt12
\font\twelvesl=cmsl12
\font\twelvesy=cmsy10 scaled \magstep1
\font\twelvei=cmmi12
\font\twelveit=cmti12
\font\twelvesc=cmcsc10 scaled \magstep1
\font\twelvebit=cmssi12
\font\twelvebbb=msym10 scaled \magstep1
\font\tenrm=cmr10
\font\tenbf=cmb10 
\font\tentt=cmtt10 
\font\tensl=cmsl10 
\font\tensy=cmsy10 
\font\teni=cmmi10 
\font\tenit=cmti10 
\font\tensc=cmcsc10 
\font\tenbit=cmssi10
\font\tenbbb=msym10
\font\ninei=cmmi9 
\font\ninerm=cmr9

\font\ninesy=cmsy9 
\font\eighti=cmmi8
\font\eightrm=cmr8

\font\eightsy=cmsy8
\font\seveni=cmmi7 
\font\sevenrm=cmr7
\font\sevensy=cmsy7
 
%
% large fonts
%

%

\newfam\bbbfam

% this macro defines a tenpoint font
\def\tenpoint{
\def\rm{\fam0\tenrm}
\textfont0=\tenrm \scriptfont0=\eightrm \scriptscriptfont0=\sevenrm
\textfont1=\teni \scriptfont1=\eighti \scriptscriptfont1=\seveni
\textfont2=\tensy \scriptfont2=\eightsy \scriptscriptfont2=\sevensy
\textfont3=\tenex \scriptfont3=\tenex \scriptscriptfont3=\tenex

\textfont\itfam=\tenit
\def\it{\fam\itfam\tenit}

\textfont\slfam=\tensl
\def\sl{\fam\slfam\tensl}

\textfont\bffam=\tenbf
\def\bf{\fam\bffam\tenbf}

\textfont\ttfam=\tentt
\def\tt{\fam\ttfam\tentt}

\def\sc{\tensc}

\def\bit{\tenbit}

\def\bbb{\fam\bbbfam\twelvebbb}
\textfont\bbbfam=\tenbbb
\scriptfont\bbbfam=\eightrm 
\scriptscriptfont\bbbfam=\sevenrm
 
\normalbaselineskip=12pt

\setbox\strutbox=\hbox{\vrule height10pt depth4pt width0pt}%
\normalbaselines\rm}

% this macro defines a twelvepoint font

\def\twelvepoint{
\def\rm{\fam0\twelverm}
\textfont0=\twelverm \scriptfont0=\ninerm \scriptscriptfont0=\sevenrm
\textfont1=\twelvei \scriptfont1=\ninei \scriptscriptfont1=\seveni
\textfont2=\twelvesy \scriptfont2=\ninesy \scriptscriptfont2=\sevensy
\textfont3=\tenex \scriptfont3=\tenex \scriptscriptfont3=\tenex

\textfont\itfam=\twelveit
\def\it{\fam\itfam\twelveit}

\textfont\slfam=\twelvesl
\def\sl{\fam\slfam\twelvesl}

\textfont\bffam=\twelvebf
\def\bf{\fam\bffam\twelvebf}

\textfont\ttfam=\twelvett
\def\tt{\fam\ttfam\twelvett}

\def\sc{\twelvesc}
\def\bit{\twelvebit}

\def\bbb{\fam\bbbfam\twelvebbb}
\textfont\bbbfam=\twelvebbb
\scriptfont\bbbfam=\ninerm 
\scriptscriptfont\bbbfam=\sevenrm
 
\normalbaselineskip=14pt

\setbox\strutbox=\hbox{\vrule height10pt depth4pt width0pt}%
\normalbaselines\rm}

% this macro defines a fourteen point font

\def\fourteenpoint{
\def\rm{\fam0\fourteenrm}
\textfont0=\fourteenrm \scriptfont0=\twelverm \scriptscriptfont0=\tenrm
\textfont1=\fourteeni \scriptfont1=\twelvei \scriptscriptfont1=\teni
\textfont2=\fourteensy \scriptfont2=\twelvesy \scriptscriptfont2=\tensy
\textfont3=\tenex \scriptfont3=\tenex \scriptscriptfont3=\tenex

\textfont\itfam=\fourteenit
\def\it{\fam\itfam\fourteenit}

\textfont\slfam=\fourteensl
\def\sl{\fam\slfam\fourteensl}

\textfont\bffam=\fourteenbf
\def\bf{\fam\bffam\fourteenbf}

\textfont\ttfam=\fourteentt
\def\tt{\fam\ttfam\fourteentt}

\def\sc{\fourteensc}
\def\bit{\fourteenbit}

\def\bbb{\fam\bbbfam\twelvebbb}
\textfont\bbbfam=\fourteenbbb
\scriptfont\bbbfam=\tenrm 
\scriptscriptfont\bbbfam=\eightrm

\normalbaselineskip=16pt

\setbox\strutbox=\hbox{\vrule height10pt depth4pt width0pt}%
\normalbaselines\rm}

% Set some default dimensions
%
\twelvepoint

\abovedisplayskip 14pt plus 3pt minus 10pt%
\belowdisplayskip 14pt plus 3pt minus 10pt%
\abovedisplayshortskip 0pt plus 3pt%
\belowdisplayshortskip 8pt plus 3pt minus 5pt%
\parskip 3pt plus 1.5pt
\hsize=6.5in
\vsize=8.9in

\fi			% end of \ifx\UsualIsLoaded

\def\ANG{Angenent}
\def\Ca{C^\ast}
\def\DB{dual billiards}
\def\DD{Douady}
\def\G{\Gamma}
\def\HTM{half-cylinder twist map}
\def\LL{{\cal L}}
\def\Lip{Lipschitz}
\def\Leb{Lebesgue}

\def\QED{  \rlap{$\sqcup$}$\sqcap$ }
\def\Tab{Tabachnikov}
\def\a{\alpha}
\def\arcot{\hbox{arcot}}
\def\bx{{\bar x}}
\def\cD{{\cal D}}
\def\cf{{\it cf.}}
\def\cl{\centerline}
\def\circle{S^1}
\def\d{\; d}
\def\eg{{\it eg.}}

\def\gp{\gamma^\prime}
\def\g{\gamma}
\def\ie{{\it i.e.}}
\def\integers{{\bbb Z}}
\def\inv{^{-1}}
\def\mG{$\Gamma$}
\def\mo{_{-1}}
\def\naturals{{\bbb N}}

\def\rationals{{\bbb Q}}
\def\ra{\rightarrow}
\def\reals{{\bbb R}}
\def\plane{{\bbb R}^2}
\def\tCa{{\tilde C}^\ast}
\def\tC{{\tilde C}}

\def\tO{\tilde\Omega}
\def\tf{\tilde{f}}

\def\tp{{\theta^\prime}}

\def\tz{{\tilde z}}
\def\t{\theta}
\def\htheta{{\hat \theta}}
\def\ux{\underline{x}}
\def\uy{\underline{y}}
\def\uz{\underline{z}}
\def\va{{\bf \alpha}}
\def\vn{{\bf n}}
\def\vu{{\bf u}}
\def\vv{{\bf v}}
\def\vw{{\bf w}}
\def\words#1{\ \hbox{#1}\ }
\def\xp{{x^\prime}}

\def\overbar{\bar}

\def\deta{{\dot \eta}}
\def\longfig#1#2#3{%
\midinsert\centerline{\psfig{figure=#1}}\medskip\centerline{%
\hbox to 4 in{\vbox{\hsize=3.9in\noindent {\bf Figure #2:} #3. } }}\endinsert}

\def\fig#1#2#3{%
\midinsert\centerline{\psfig{figure=#1}}\medskip\centerline{%
{\bf Figure #2:} #3. }\endinsert }

\def\sgn{\hbox{sgn}}

\input psfig

\centerline{\bf Dual Billiards, Twist Maps and Impact Oscillators}
\bigskip
\cl{Philip Boyland\footnote{ }{ The author was
partially supported by NSF grant \# 431-4591-A}}
\cl{Institute for Mathematical Sciences}
\cl{SUNY at Stony Brook}
\cl{Stony Brook, NY 11794}
\cl{Internet: boyland@math.sunysb.edu}
\bigskip
{\narrower

\noindent{\bf Abstract.} 
In this paper techniques of twist map theory are applied
to the annulus maps arising from dual billiards 
on a strictly convex closed curve $\G$ in the plane.
It is shown that there do not exist   invariant circles
near \mG\ when there is a point on \mG\
where the radius of curvature vanishes or is discontinuous.
In addition,  when the radius of curvature is not $C^1$ there are
examples with orbits that converge to a point of \mG. 
If the derivative of the radius of curvature is bounded,
such orbits cannot exist.
The final section of the paper concerns an impact
oscillator whose dynamics are the same as a  \DB\ map.
The appendix is  a remark on the connection of  the
inverse problems for invariant circles in billiards and dual billiards.

}

\bigskip
\noindent{\bf \S 0 Introduction}

Dual billiards is a dynamical system defined
on the exterior of an oriented convex closed curve \mG\
 in the plane. If $z$ is a
point in the unbounded component of $\plane - \G$, then
 its image under the \DB\ map $\Phi$ is the reflection
about the point of tangency  in the  oriented 
supporting line to \mG\ (see Figure 0.1).
It is clear that $\Phi$ is an area-preserving, and if \mG\ is
strictly convex, it is homeomorphism. 
Further, $\Phi$ is affine invariant, \ie\
transforming the convex curve by an affine map
simply transforms the entire homeomorphism.
If the initial curve is an ellipse, then its exterior 
is foliated by concentric homothetic ellipses that are $\Phi$-invariant.

\fig{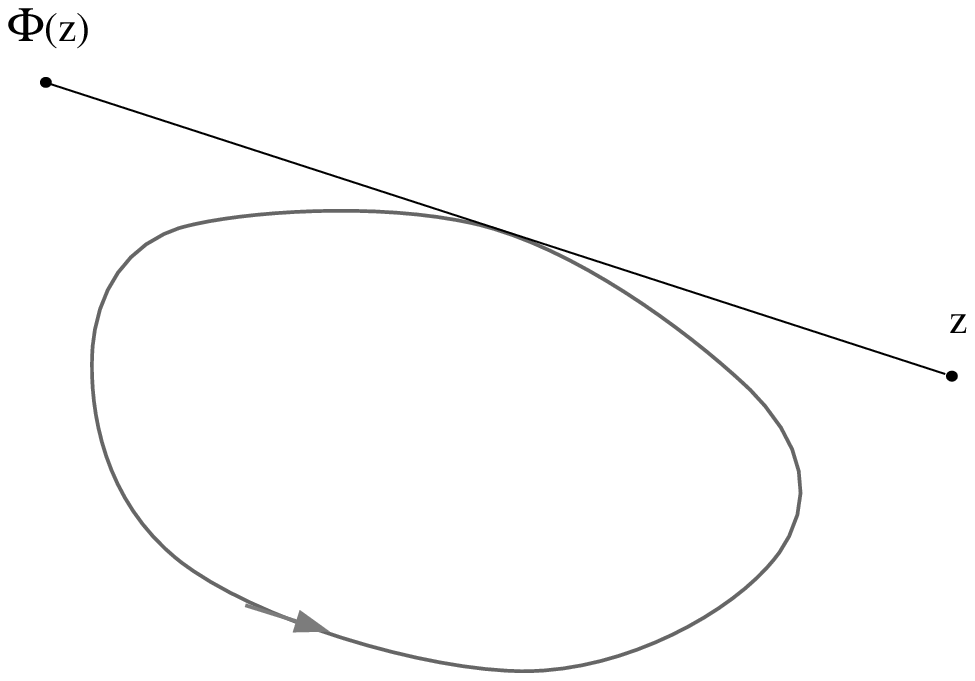,height=.30\hsize}{0.1}{The dual billiards map}

The invention of \DB\ is credited to B. H. Newman
in [Ms2]. The nomenclature 
outer billiards or cobilliards is also used.
The name \DB\ is apt  because this dynamical
system is in certain ways the dual of the usual
billiards. Indeed, many of
the theorems given here are analogs of theorems
for billiards. There are, however, fundamental differences between the
systems. Perhaps the most important is
that the phase space for the \DB\ map is non-compact 
and has infinite area. 
An obvious obstruction to any explicit duality between 
the two systems is that the
billiards map defined using an
ellipse has period-two points of elliptic and saddle type.

One of the goals of this paper is to show
how \DB\ fits into recent developments
in the  theory of area-preserving twist maps
of the annulus due mainly to Mather.
Once it is known that
a dynamical system falls into this class there are a wealth
of results available which yield a great deal of information about the dynamics of the map.
There are technical differences, however, between the
twist maps that arise in  \DB\ and the usual  classes of twist maps that 
are studied in the literature.

Section 1  begins by  specifying a class of twist maps  called
\HTM s that  includes those from \DB. 
The rest of \S 1 is devoted to describing some standard tools and results of
twist map theory and remarking on how these
can be adapted to the new  class.  We only describe that part of 
twist map theory that is relevant. For a broader perspective the 
reader is referred to the various papers of Mather (summarized in
[M5]), as well as [Bg], [Ms3], [Ms4],  chapter 1 of [H] and chapter 10 of [MH].
There are surveys with a more physical point of view in 
[Mc] and [Me].

Section 2 give necessary background material
on convex curves and in \S 3 it is confirmed that,
in the coordinates introduced in \S 2,	 the 
\DB\ map is indeed a \HTM.
The next section deals with topological 
circles that are invariant under the \DB\ map.  
The question of the existence of invariant circles
is central in twist map theory. An invariant circle
provides a barrier through which orbits cannot pass,
and thus the existence of invariant circles is intimately connected
with stability questions. Section 4 begins with two results of
R. \DD\ that give existence of invariant circles
near \mG\ and near infinity. 
These results require a certain smoothness as they ultimately
depend on KAM theory. 

The inverse problem for invariant circles
is considered next: given a convex curve
$\G_1$ is there a $\G$ whose \DB\ map has $\G_1$ as
an invariant circle.
The solution to this problem is given by 
a   geometric construction called the area envelope. The analog
of this construction for billiards 
is discussed in the appendix where it is pointed out that
both the billiard and dual billiard constructions are  a consequence of the
constancy of an area function (Mather's big
``H'') that can be associated with an invariant circle
of a twist map (\cf\ Proposition 1.3).

The next results in \S 4 concern the   non-existence of invariant circles 
near \mG\ when there is a point on \mG\
where the radius of curvature vanishes or is discontinuous.
 In addition,  when the radius of curvature is not $C^1$ there may also exist
orbits that converge to a point of \mG. These
so-called crash solutions are shown to exist for
certain examples in \S 5. It is also shown that if
the derivative of the radius of curvature is bounded,
no such orbits exist.

There are no known examples of strictly convex curves
for which the dual billiards map lacks invariant circles near infinity. 
This problem (raised by Moser) remains one of the most important outstanding 
problems in the theory of \DB. 

\fig{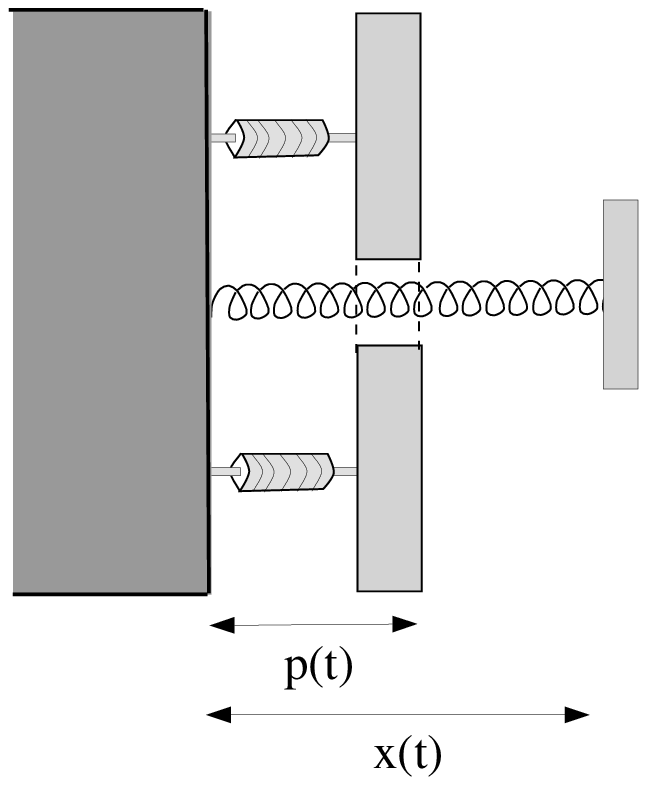,height=.30\hsize}{0.2}{The impact oscillator}

The final section of the paper concerns an impact
oscillator whose dynamics are the same
as a  \DB\ map. Consider a wall that is
periodically moving with position given by $p(t)$ where
$p(t)$ satisfies
$\ddot p + p = \rho(t)$ for some $2 \pi$-periodic forcing function
$\rho$ with $\int \exp(i t) \rho(t) = 0$.
A simple harmonic oscillator with position specified by
$x(t)$ satisfying $\ddot x + x = 0$ 
is assumed to have perfectly elastic collisions
with the wall. Figure 0.2 shows one possible (idealized) realization
of the system. On the left is a fixed wall whose position is set at
the origin. A pair of pistons periodically drive a second wall that is
to the right of the fixed wall. A flat plate is connected to the fixed wall
by a spring that passes through a  hole in the moving wall.
We monitor the motions of the plate as it moves away, collides with
the moving wall, moves away, etc. The return map to collisions
with the moving wall in appropriate coordinates
turns out to be  identical to \DB\ on
a convex curve $\G$ whose radius of curvature function $\rho$
is the same as the forcing function for the wall.

As a consequence of the equivalence of the
two dynamical systems one can use  results
about \DB\ to provide information about the oscillator. 
For example,
as a consequence of the Birkhoff-Mather stability theorem
(see Remark 4.3.2) when there do not exist $\Phi$-invariant circles
near the curve $\G$, there must exist orbits that converge
to $\G$ in forward and backward time. Thus as a consequence of 
Theorem 4.3, if there
is an instant at which the forcing function is zero or
discontinuous,
then there are solutions for the impact oscillator that converge to the wall.
In addition, it follows from Theorem 5.1 that there are impact 
oscillators that have  crash solutions that
converge to the wall in finite time. These solutions do not
exist if the forcing function is $C^1$ (Theorem 5.2).
 
The situation near infinity is perhaps of more physical interest.
As a consequence of Theorem 4.1, for sufficiently smooth convex
curves the \DB\ map has invariant circles near infinity.
Thus for sufficiently smooth forcing functions the 
impact oscillator has no unbounded solutions.
Because of the lack of unbounded solutions 
these oscillators are called {\it stable}.
There have been numerous other impact oscillators studied in the
literature. Their stability depends on the details
of the oscillator, \eg\ on the restoring force of the
particle to the wall, the relation of
the frequency of the forcing to that
of the free oscillator, etc.
Again using the Birkhoff-Mather stability theorem, the
non-existence of invariant circles
near infinity for \DB\ on some convex curve
would imply unbounded solutions for the corresponding impact oscillator.
This gives a physical motivation for the question of Moser noted above.

This paper has several purposes, and it
contains material of different types.
A first goal is to  survey known results from
the  point of view adapted here. For this reason there is an
overlap in basic material with other
papers that have appeared on \DB, for example,
[GK], [T1], [T2], and [T3].
A second goal is to give  a mathematical framework for future work on \DB. 
Thus there are somewhat technical sections 
(\eg\ \S 1.1) that
could be  skipped by the reader with more physical interests. 
Despite these rigorous goals we also want to preserve the essentially
simple and intuitive geometric character of the
problem. Thus there also is descriptive
material. Much of this material also serves  to illustrate
geometric heuristics associated with twist map theory.

\bigskip
\noindent{\bf Acknowledgments:} Most of the results in this paper were
obtained while the author was a Postdoctoral Research Fellow
at the Mathematical Sciences Research Institute in Berkeley,
1988-89. He would like to thank them (belatedly) for there support.
Thanks also to Danny Goroff, Dick Hall, and Robert MacKay
for useful and stimulating conversations.

\bigskip

\noindent{\bf \S 1 Half-cylinder twist maps}

In this section we introduce a class of area-preserving
monotone twist maps called half-cylinder twist maps,
and point out that the standard results from twist map
theory apply to this class.
We will see in \S 3 that a \DB\ map
in the appropriate coordinates is a half-cylinder twist map.

\medskip
\noindent{\bf \S 1.1 Definition of half-cylinder monotone twist maps.}
Let the one point compactification of the ray $[0, \infty)$ 
be denoted $[0, \infty]$ with the added point being labeled $\infty$.
Let $C$ be the open  half-cylinder $\circle\times (0,\infty)$
and $\Ca$ be $\circle\times [0,\infty]$, with the circle
$\circle=\reals/(2 \pi \integers)$. The space $\Ca$ is homeomorphic
to the compact annulus, but we maintain the infinite measure.
The universal covers of these spaces are denoted by $\tC$ and
$\tCa$, respectively. A tilde will always indicate the lift of
a point, function, set, etc. to the universal cover.
 The map $\pi_1$ in various contexts 
means the projection of a product space onto its first component.
We will use the coordinates $(\t, \gamma)$ for $C$ and $\Ca$  and 
$(x, \gamma)$ for $\tC$ and $\tCa$.
If $g$ is a real valued function,
$g_i$ denotes the derivative with respect to the $i^{th}$ variable.

The definition of a half-cylinder monotone twist map requires two pieces
of data. The first is a pair of functions $b_0, b_\infty
:\reals\ra\reals$ called the {\de boundary maps}.
It is required that $b_0$ be nondecreasing and continuous 
from the right, $b_\infty$ be increasing and continuous,
$b_0 < b_\infty$, and  $b_i(x+1) = b_i(x) + 1$
for $i = 0, \infty$.
Let $\cD = \{ (x, \xp)\in \reals^2 : b_0(x) < \xp < b_\infty(x)\}$.
For $i = 0,\infty$, $\Lambda_i\subset \reals^2$ denotes the graph
of $b_i$ as a function of $x$. Note that $\Lambda_0\cup\Lambda_\infty$
is contained in the frontier of $\cD$. 

The second piece of data needed to define a half-cylinder twist
map is a $C^1$-function $h:\cD\ra \reals$ called
a {\de generating function} that satisfies:

(G1) (periodicity) $h(x + 2\pi, \xp + 2\pi) = h(x, \xp)$

(G2) (conditions on the derivatives) 
 $h_1 < 0$ and $h_2 > 0$, and for each $\epsilon > 0$, the derivative $Dh$ is 
Lipschitz on each 
$\cD_\epsilon = \{ (x, \xp)\in \reals^2 : b_0(x) < \xp < b_\infty(x)- \epsilon\}$.
Further, the mixed partial derivative $h_{12} = h_{21}$ exists, 
is continuous and satisfies $h_{12} < 0$.

(G3) (limit behavior near $\Lambda_0$) For $i = 1,2$,
$$\lim_{(x, \xp)\ra \Lambda_0} h_i(x, \xp) = 0 = 
\lim_{(x, \xp)\ra \Lambda_0} h(x, \xp).$$

(G4) (limit behavior near $\Lambda_\infty$) The first derivatives must satisfy
$$\lim_{\xp\ra b_\infty(x)} h_2(x, \xp) = \infty = \lim_{x\ra b_\infty^{-1}(x)}
h_1(x, \xp)$$  
and for each $x$,
$$
h_2(u, v) + h_1(v, w) \ra 0 \words{as} 
(u,v,w)\ra (b_\infty\inv(x), x, b_\infty(x)).
$$

\medskip

{\bf Definition:} {\it  A map $f:\Ca\ra \Ca$ is called a {\de half-cylinder
monotone twist map} if it has a lift $\tf:\tCa\ra\tCa$ that 
satisfies
$$\eqalign{\g &= - h_1(x, \xp)\cr
\gp &= h_2(x, \xp)\cr}\eqno(1.1)$$
where $(\xp, \gp ) = \tf(x,\g )$ for $(x, \gamma)\in \tC$,
and $h$ satisfies properties (G1) -- (G4)
above with $b_i$ equal to $\tf$ restricted to $\reals\times \{i\}$.
}
\medskip

As a consequence of the properties (G1) -- (G4), when restricted
to $C$, a half-cylinder
monotone twist map $f$ will always be an area-preserving  homeomorphism,
and both $f$ and $f^{-1}$ will be locally Lipschitz.
Conversely, given $b_i$ and $h$ satisfying (G1) -- (G4), 
(1.1) implicitly defines a half-cylinder monotone twist map
$f$.
Using (1.1), $D\tf (x,\g)$ is given by
$${1\over h_{12}}
\left(\matrix{ -h_{11} & -1\cr
  h_{12}^2 - h_{11} h_{22}  & -h_{22}\cr} \right) \eqno(1.2)
$$
where each of the derivatives $h_{ij}$ is evaluated at
$(x,\xp)$ with $\xp = \pi_1(\tf(x, \g))$.

{\bf Remarks}

{\bf 1.1.1} The name ``twist map'' is given because the 
image of a vertical arc is
always the graph of a function defined on a subset
of the circle. From (1.2) one sees that
this is a consequence of $h_{12} <0$.
The generating function $h$ can be given a geometric
interpretation as in Figure 1.1 ([K]).
The area of the region bounded by the $x$ axis,
the vertical arc above $\xp$, and
the image under $\tf$ of the vertical arc above $x$ is
equal to $h(x, \xp)$.

\longfig{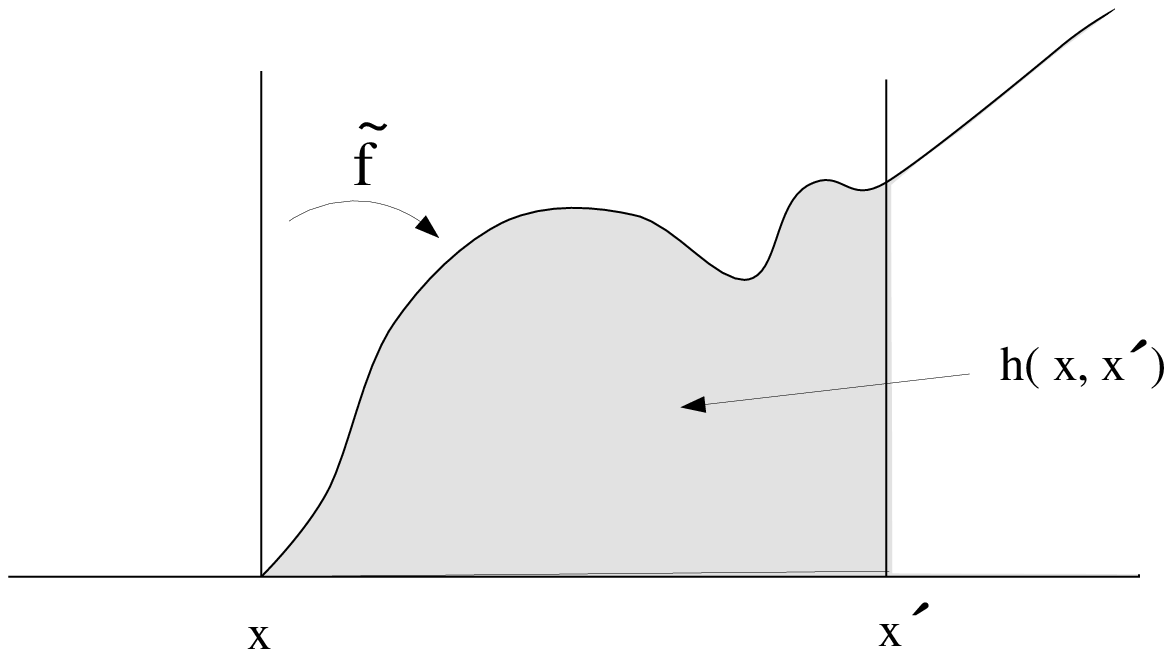,height=.30\hsize}{1.1}{The geometric meaning of the generating function}

{\bf 1.1.2} There are a few differences
between half-cylinder twist maps and the classes of twist maps that
occur in the literature.
These differences are mainly ``technical'' in the sense that
standard arguments work for the class with straightforward alterations.
The consideration of these technical differences are forced on us by
the twist maps that come from \DB.

A twist map is usually defined either on the compact annulus 
or on the infinite cylinder.
In the later case one requires infinite wrapping, 
\ie\ the image of a vertical line wraps infinitely around
the cylinder. In contrast, \HTM\ are defined
on the half-infinite cylinder. They extend to
a perhaps non-continuous map on the lower boundary, and
can be extended to a homeomorphism on
an upper circular boundary that is added at
infinity. Thus they do not have infinite wrapping;
they can be thought of as being
defined on a compact annulus with infinite area.

If a twist map $f$ is differentiable at a point $z_0$
the {\de positive twist} at $z_0$, denoted $twist^+(z_0)$,
 is the angle (measured clockwise) from a vertical vector based at 
$f(z_0)$ to the image under $Df$ 
of a vertical vector based at $z_0$.
The lift of a  \HTM\ $\tf$ may not be differentiable, but the image of the
vertical line through $\tilde z_0 = (x_0, \g_0)$ is always the 
graph of the Lipschitz
function of $\xp$, $h_2(x_0, \xp)$. In this case 
we let $twist^+(z_0)= \arcot(K)$ where $K$ is the best \Lip\
constant for this function at $\xp_0 = \pi_1(\tf(z_0))$.
The {\de negative twist} at $z_0$, denoted 
$twist^-(z_0)$, is the twist at $z_0$ of $S\circ \tf$ 
where $S(x, \g) = (-x, \g)$.
The {\de twist} at $z_0$ is $twist(z_0) = \min\{ twist^+(z_0),
twist_-(z_0)\}$. Note that if $h$ is twice differentiable,
$twist(z_0)= \min\{\arcot(h_{11}(\hat x_0, x_0),\arcot(h_{22}(x_0, \xp_0)\}$,
where $\hat x_0 = \pi_1(\tf\inv(z_0))$.

In twist map theory it is usually assumed that the
twist is bounded away from zero. For \HTM s the
condition given in (G2) requiring
that  $Dh$ is \Lip\ in $\cD_\epsilon$
insures that the twist is bounded away from
zero in neighborhoods of the lower boundary of $C$. 
However,   $f$ itself may only be
only locally \Lip\ as its \Lip\ constant
can go to infinity as one approaches the lower boundary.
In a \HTM\ the twist will not be bounded away from
zero as one approaches infinity
as the map extends to the circle at infinity.
\medskip

\noindent{\bf \S 1.2 The variational formulation and dynamics of twist maps.}
Note that (1.1) says that $h$ is a generating function of a canonical
transformation in the sense of classical mechanics (\eg\ [MH], page 47).
The function $h$ plays another role 
and is  sometimes called the {\de action} as it can be used in a discrete
variational formulation. What follows is the rudiments of this
theory adapted for \HTM s. For a full treatment the reader is
referred to [Bg].

A {\de configuration} is an element $\ux = (x_i)_{i\in \integers}$
of the set $\reals^\integers$ of bi-infinite sequences of real numbers
(the nomenclature ``configuration'' comes from solid state physics).
Given a generating function $h$ and 
boundary maps $b_i$ as above, let $X= \{ \ux \in \reals^\integers :
(x_i, x_{i+1})\in \cD \hbox{ for all } i \}$.
Elements of $X$ are called {\de allowable configurations}.
The generating function $h$ can be extended to finite segments
of allowable configurations  via
$$h(x_j, \dots\; , x_k) = \sum_{i=j}^{k-1} h(x_i, x_{i+1}).$$
An allowable configuration $\ux$ is called {\de stationary} if
for all pairs $(j, k)$ with $j < k$,
$h(x_j, \dots\; , x_k)$ is stationary with respect to variations with
fixed endpoints. The configuration is called
{\de minimizing} if
for all pairs $(j, k)$ with $j < k$,
 $h(x_j, \dots\; , x_k)$ is minimal with respect to variations with
fixed endpoints. Minimizing configurations are
always stationary. A configuration is stationary  if and only if
$$h_2(x_{i-1}, x_i) + h_1(x_i, x_{i+1})= 0\eqno(1.3)$$
for all $i$,
or as a consequence of (1.1), precisely when $\ux$ is a sequence
of $x$ coordinates of an orbit of the map $\tf$, \ie\  for all $i$,
$x_i = \pi_1(\tf(x_0, \g_0))$ for some initial $(x_0, \g_0)\in \tC$.

A configuration is {\de monotone} if for each fixed $(m,n)\in \integers^2$
either $x_{i+m} > x_i + n$ for all $i$, or else,
$x_{i+m} < x_i + n$ for all $i$.
A basic feature of the theory is that minimizing configurations
are always monotone. An orbit $(x_i, \g_i)$ of $\tf$ is
{\de monotone} if the sequence $(x_i)$ is. An orbit 
$(\t_i, \g_i)$ of $f$ is {\de monotone} if it has a lift to
$\tC$ that is monotone. Roughly speaking, a monotone orbit is one
on which $f$ preserves the radial order. An invariant set is called
monotone if each of its orbits is.

Given the lift of a half-cylinder twist map $\tf$
the {\de rotation number} of a point $\tz\in \tCa$ is
$$rot(\tz, \tf) := \lim_{n\ra\infty}
{\pi_1(\tf^n(\tz)) - \pi_1(\tz)\over 2 \pi n}$$
if the limit exists. The rotation number of a
point $z\in \Ca$ is the rotation number of one (and thus all) of
its lifts. Note that this is
defined only up to a choice of lift of $f$, \ie\  up
to integer translation. The {\de rotation set}
of $\tf$ is $rot(\tf) := \{ rot(\tz, \tf) : \tz\in \tCa\}$. 
The rotation number of
a point under a circle map is defined similarly, only in this case,
if $g:\circle\ra\circle$ is nondecreasing (but not necessarily
continuous), then the rotation number exists and is the same for
all points in the circle. This common number
is called the {\de rotation number} of the map, and is denoted
$rot(g)$. Another basic feature of twist map theory is that
monotone orbits always have a rotation number. This is basically because
they behave like orbits of 
circle maps. Note that any nondecreasing $G:\reals\ra \reals$ that
satisfies $G(x + 2 \pi) = G(x) + 2 \pi$ may also be assigned a rotation
number that will be denoted $rot(G)$.

A Denjoy minimal set for a half-cylinder twist map
is a compact, invariant set $X$ so that $f$ restricted
to $X$ is conjugate to the exceptional
minimal set in a Denjoy counterexample on the circle 
(see \eg\ [Dv], page 108).
The phrase ``invariant circle'' in this paper
always means a {\it homotopically nontrivial}
circle that is invariant under the map. Since $f$ restricted
to an invariant circle, $\Omega$, is a circle homeomorphism, it has a rotation
number that will be denoted $rot(\Omega)$.

\medskip
{\bf Theorem 1.1:} {\it If $f$ is a half cylinder twist map with
boundary functions  $b_i$, then:

(a) (Aubry-Mather) The rotation set of $f$ is 
$ rot(f) = [rot(b_0), rot(b_\infty)]$.
For each $p/q\in rot(f)$ with $p$ and $q$ relatively
prime, there is  a monotone periodic
orbit with rotation number $p/q$ and period $q$. For each
$\omega\in rot(f) - \rationals$,  there is an invariant circle
or  Denjoy minimal
set $X_\omega$ that is monotone and each $x\in X_\omega$ has rotation
number $\omega$.

(b) (Birkhoff) If $\Omega \subset C$ is a $f$-invariant
circle, then $\Omega$ is the graph of a Lipschitz function
$\circle\ra (0,\infty)$. Further, if $\cal C$ denotes
the set of invariant circles (including $\circle\times \{ 0 \}$)
 with the Hausdorff topology,
then $\cal C$ is closed and the map \ $rot: \cal C \ra \reals$
is continuous and monotonic in the sense that 
$\Omega_2$ being contained in the open annulus
bounded by  $\Omega_1$ and $\circle\times \{0\}$ implies
$rot(\Omega_1) > rot(\Omega_2) > rot(b_0)$.

(c)  (Mather) If $f$ has a invariant circle $\Omega\subset C$,
then for almost every $(x_0,\g_0)\in \Omega$, 
$$h_{11}(x_0, x_1) + h_{22}(x_{-1}, x_0) >0$$
where $(x_{i}, \g_{i}) = f^{i}(x_0, \g_0)$.
}
\medskip
{\bf Proof:} Part (a) can be proven using variational methods
in which case the monotone orbits are obtained as
minimizing orbits (\eg\ [Bg]). It may also be proven
using topological methods
(\eg\ [MH]), or using monotone recursion
relations ([A]). In part (b), the statement that invariant
circles are Lipschitz graphs can be  proved as in [M1] and [M2],  or 
chapter 1 of [H]. The statement that the set of invariant circles
is closed in the Hausdorff topology requires the observation
(explained in Remark 1.2.2 below) that condition (G2) implies a uniform
bound on
the \Lip\ constant of functions $u:\circle\ra (0, \infty)$
whose graphs are invariant circles in $\circle\times [0, K]$.
The monotonicity of the rotation numbers on the circles is a standard
consequence of the twist hypothesis.

To prove (c),  let $\alpha$ denote the homeomorphism of
the real  line that is first component of $\tf$ restricted to 
the lift of the invariant circle, \ie\ $\alpha(x) = \pi_1(\tf(x, u(x)))$.
 Since $\tf$ is locally \Lip\ and $u$ is Lipschitz, $\alpha$ and 
$\alpha\inv$ are \Lip, and thus their derivatives exist almost everywhere. 
Equation (1.2) yields that 
$$h_1(x, \alpha(x)) + h_2(\alpha\inv(x), x)  \equiv 0. \eqno(1.4)$$
 Part (c) then follows by differentiating and using  $h_{12}< 0$.\QED
\medskip

{\bf Remarks:}

{\bf 1.2.1} A geometric interpretation of the condition
in  Theorem 1.1(c) can be  given in terms of ``tumbling tangents''
(\cf\ [Mc] and [MP]). If $h$ is twice differentiable, using (1.2) one sees
that  the slope of the image of the vertical tangent vector
at the point $(x, \g)$ is given by $h_{22}(x, \xp)$, 
where $\xp = \pi_1(\tf(x, \g))$.
Similarly, under $\tf\inv$, the slope of the
image of the vertical tangent vector
at the point $(\xp, \g^\prime)$ is $-h_{11}(x, \xp)$.

Now assume there is  an invariant circle that passes 
through the three points $(x_i, \g_i)$
for $i = -1, 0, 1$ where $(x_{i+1}, \g_{i+1}) = \tf(x_i, \g_i)$.
Theorem 1.1(b) says that the circle is the graph of a 
Lipschitz function, and thus the unit tangent vector in the direction of
the circle exists almost everywhere; let us 
assume that it exists at our triple. Since the circle is invariant,
the collection of these circle directions is invariant under
the induced action of $D\tf$ on the unit tangent bundle.
Now since $\tf$ is a twist map, the image of a unit vertical
vector always points in the positive $x$ direction.
Because the three points are on an invariant circle and the circle
directions are invariant, the second iterate of the 
vertical vector at $(x\mo, \g\mo)$ cannot rotate
past the downward vertical.

\longfig{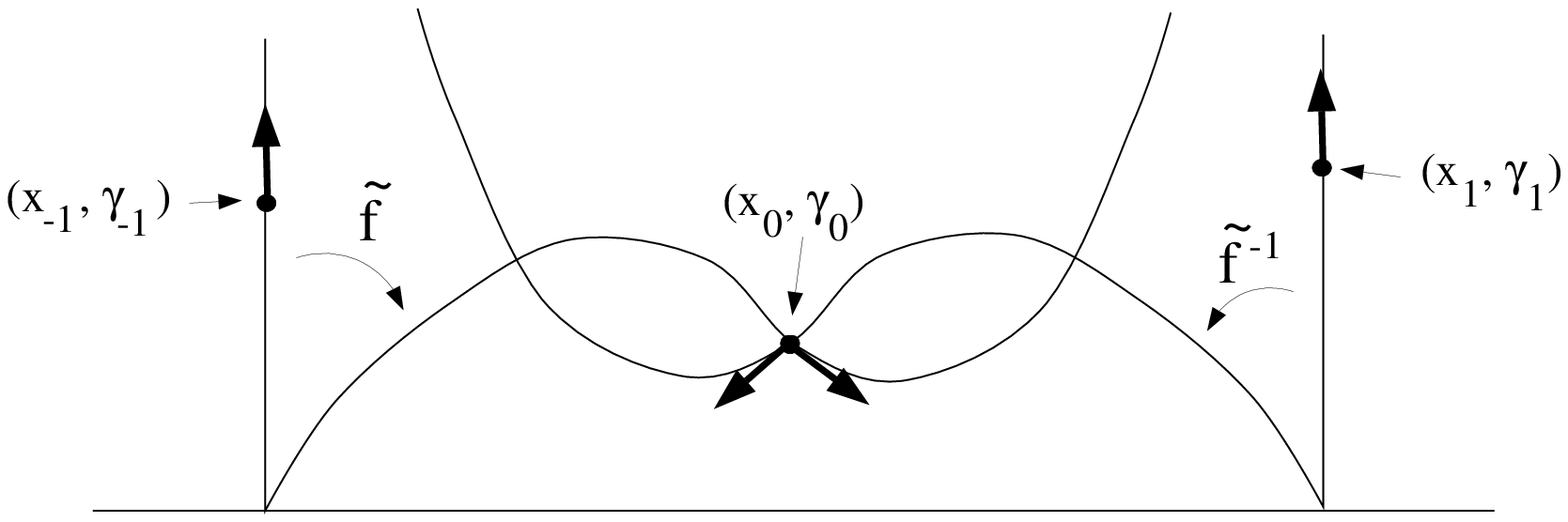,height=.25\hsize}{1.2}{The tumbling tangents criterion for the non-existence
of invariant circles}

Rather than apply these considerations to the second iterate,
it is usually more convenient to  iterate forward from $(x\mo, \g\mo)$
and backwards from $(x_1, \g_1)$. Since the slope
of the forward image of a vertical vector is $h_{11}$ 
and a backward image is $-h_{22}$, the ``no tumbling tangents''
condition is  exactly that given in Theorem 1.1(c).  Figure 1.2
shows an example in which there can be no invariant circle
containing the three points.

By using more iterates, these considerations can be refined 
to obtain the converse KAM theory of MacKay and Percival
([MP]).
The macroscopic version of tumbling tangents is the existence
of non-monotone periodic orbits as in [B+H].

{\bf 1.2.2} The regularity conditions
imposed on \HTM s were chosen so that the standard proofs of
Theorem 1.1(b) and (c) would go through with little change. 
In order to differentiate
(1.4) almost everywhere, we need that $f$ is locally
\Lip.  For  this it is enough that
$Dh$ be locally \Lip. We also need $u$ to be  \Lip. 
For this the proofs ([M1], [M2], [H]) require only $f$ is a locally
 \Lip\ twist map
whose twist is bounded away from zero in a neighborhood of
the circle. Again this requires only that 
$Dh$ be locally \Lip. 

However, to get  a uniform bound in neighborhoods
of the lower boundary for the 
\Lip\ constants of the  functions $u:\circle\ra (0,\infty)$
that define invariant circles, we need 
a uniform bound on the twist (as defined in Remark 1.1.2).
This requires that $h_{11}$ and $h_{22}$ be bounded and hence
the requirement in (G1) that
$Dh$ is \Lip\ in $\cD_\epsilon$.

One can geometrically see the dependence of the \Lip\ constants of the
functions $u$ by thinking 
of the twist as defining a cone about 
the vertical in the tangent bundle at each point $x$. 
By definition, the image under $D\tf$ of a vertical 
vector based at $\tf\inv(x)$  misses this cone as does the image
of the vertical vector based at $\tf(x)$ under $D(\tf\inv)$.
This implies that
vectors in this cone are rotated beyond the vertical in either
forward or backward iteration. Thus applying the considerations
of the previous remark, no tangent to an invariant circle
can lie in this cone. 

\medskip
\noindent{\bf \S 1.3 Monotone recursion relations.}
The point of view of Angenent  will be used to 
construct examples in \S 5.1. 
The basic idea in [A] is to use (1.3) directly to define a recursion
relation that must be satisfied by sequences of $x$ coordinates
of orbits.

We will state these ideas just in the limited context needed here.
Assume that  $h$ and $b_i$ satisfy the properties (G1) -- (G4) above
and let 
$$
\Delta(x\mo, x_0, x_1) =-\bigl( h_1(x_0, x_1) + h_2 (x_{-1}, x_0)\bigr).
\eqno(1.5)
$$
The properties of $\Delta$ 
that cause it to behave like the monotone recursion relations in
[A] are: (1)   (monotonicity) $\Delta$ is a nondecreasing
function of $x\mo$ and $x_1$, (2) (periodicity) 
$\Delta(x\mo, x_0, x_1) =  \Delta(x\mo + 2 \pi, x_0 + 2 \pi, x_1 + 2 \pi)$,
 and (3) (coerciveness) which is given by the condition
in the first part of (G4).
The coerciveness condition differs somewhat from that adopted in [A], 
but the proof of Theorem 1.2 below can be constructed with little difficulty.

An allowable configuration $\ux$ is called a {\de solution} to the monotone
recursion relation given by $\Delta$ if 
$\Delta(x_{i-1}, x_i, x_{i+1}) = 0$ for all $i$. Note that
a solution always gives the sequence of $x$ coordinates of
an orbit of the lift $\tf$ of the \HTM\ defined by 
the generating function $h$. 
The configuration $\ux$ is a {\de subsolution} if 
$\Delta(x_{i-1}, x_i, x_{i+1})\allowbreak \ge 0$ for all $i$
and a {\de supersolution} if
$\Delta(x_{i-1}, x_i, x_{i+1}) \le 0$ for all $i$.
The standard partial order on configurations is given by
$\ux \le \uy$ if and only if $x_i \le y_i$ for all $i$.
The main result we need is:

\medskip
{\bf Theorem 1.2:} (Angenent) {\it Let $\Delta$ be
defined by (1.5) for a generating function $h$ associated
with a \HTM. If $\ux$ and $\uy$ are a 
subsolution and a supersolution, respectively, of
the monotone recursion relation generated by $\Delta$,
and further, $\ux \le \uy$, then there exists a solution
$\uz$ with $\ux\; \le\; \uz\; \le\; \uy $.}
\medskip

{\bf \S 1.4 The area function of an invariant circle.}
The next result concerns the constancy of a certain
 area that can be 
associated with an invariant circle. This quantity,
sometimes called ``Mather's big H'', has
a simple geometric interpretation in billiards and 
dual billiards (see the appendix and Theorem 4.2).

Let $\Omega$ be a (not necessarily invariant) circle in
$C$ that is the graph of
$u:\circle\ra (0,\infty)$. We will also use $u$ for the 
function $\reals\ra (0,\infty)$ that has  $\tilde \Omega\subset
\tC$ as its graph.  Given the lift of a  half-cylinder twist
map $f$, we can define two functions associated
with the displacement of $\tilde\Omega$ under $\tf$.  
Let  $\alpha:\reals\ra\reals$ be defined as $\alpha(x)
= \pi_1(\tf(x, u(x)))$, or equivalently, $\alpha(x)$ is the unique
solution to $-h_1(x,\alpha(x)) = u(x)$. And let
$\beta:\reals\ra\reals$ be defined by $\beta(x) =
\pi_1(\tf(x, \sigma(x)))$ where $\sigma(x) = \sup\{\gamma:
(x, \gamma) \in \tf\inv(\tilde \Omega)\}$. Equivalently,
$\beta(x) = \sup\{x^\prime : h_2(x, x^\prime) = u(x^\prime)\}$
(see Figure 1.3).

\fig{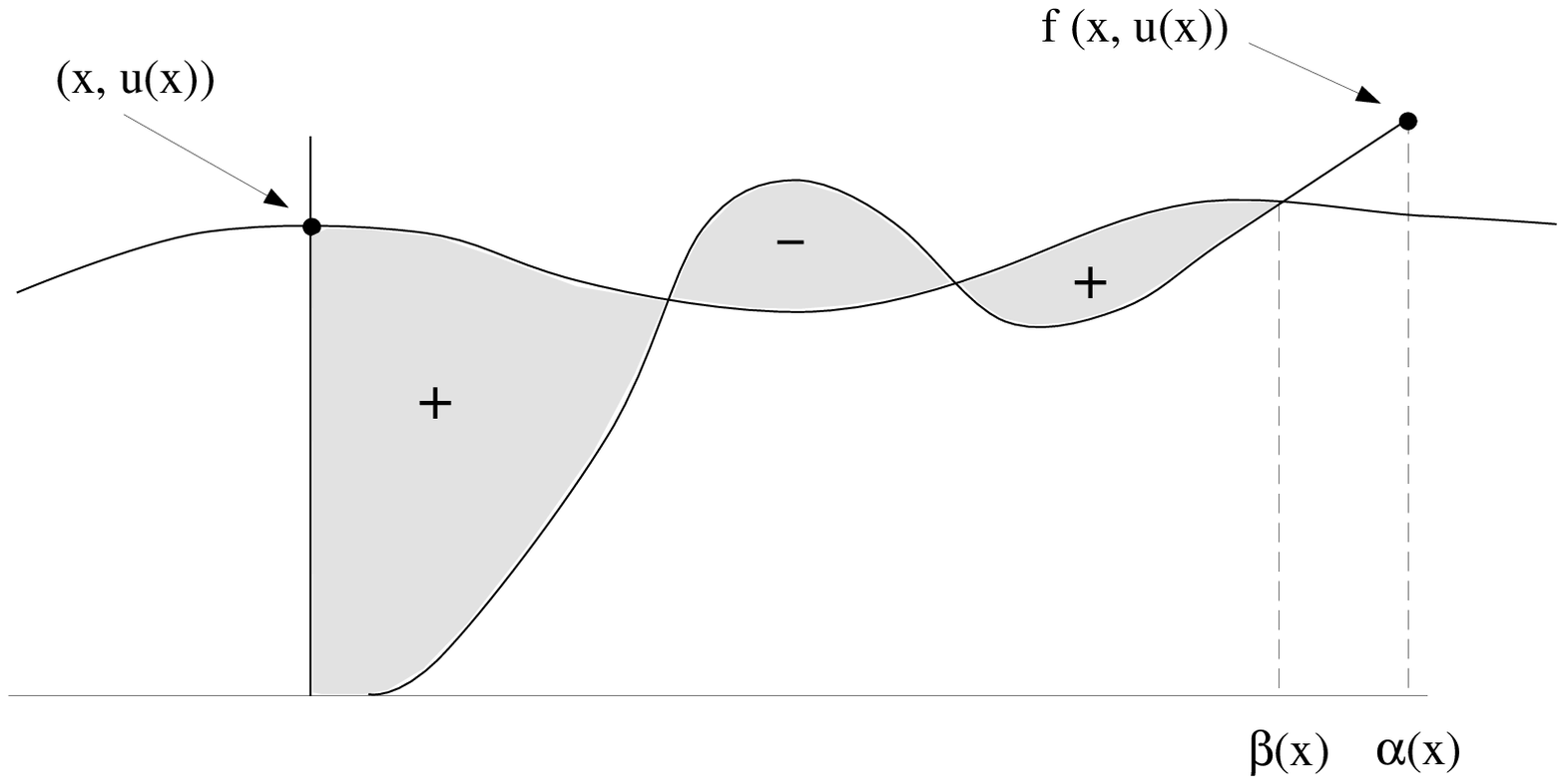,height=.30\hsize}{1.3}{The definitions of $\alpha$, $\beta$, and the area function}

From the definitions, $\alpha$ is continuous and $\beta \le \alpha$.
Since $h_{12} < 0$, $\beta$ is strictly increasing.  Note that
if $\eta > \beta(x)$, then since $\lim_{\xp\ra b_\infty(x)} = \infty$, 
we have $h_2(x, \eta) > u(\eta)$. 
Also note that $\tilde \Omega$ is $\tf$-invariant if and
only if $\alpha = \beta$.

If the generating function associated with $f$ is $h$,
define the {\de area function} of
$\Omega$ as
$$A(x) := \int_x^{\beta(x)} u(\eta) \; d\eta - h\bigl(x, \beta(x)\bigr).$$
Thus A(x) measures the signed area between the $\Omega$ and the
image of the vertical arc above $x$ (see Figure 1.3 again).
Mather proves in [M3] that the area function of an invariant circle
is always constant.

\medskip
{\bf Proposition 1.3:} (Mather) {\it Given a half-cylinder twist
map $f$,  the graph $\Omega$ 
of a continuous function $u:\circle\ra (0,\infty)$
is an $f$-invariant circle if and only if its 
area function is constant.}
\medskip

{\bf Proof:} First assume that $\Omega$ is $f$-invariant.
By Theorem 1.1(b), $u$ is a Lipschitz function and 
since $f$ is locally  Lipschitz, 
$\alpha=\beta$ is Lipschitz and so its derivative exists
almost everywhere. Differentiating $A(x)$ using its definition
and then applying (1.2) one gets that  $dA/dx = 0$ almost everywhere,
and since $A$ is continuous, it is constant.

Now conversely, assume that $A$ is constant and proceed by contradiction.
 If $\alpha \not = \beta$, then since $\alpha$ is continuous and
$\beta\le\alpha$, there exists $x$ and $\bx$ with $x < \bx$ and 
$\alpha(\eta) > \beta(\bx)$ for all
$\eta \in (x, \bx)$. Thus for these $\eta$,
$h_1(\eta, \beta(\bx)) > h_1(\eta, \alpha(\eta)) =
- u(\eta)$. Note also that above we showed that
$h_2(x, \eta) > u(\eta)$ for all $\eta > \beta(x)$.
Now by assumption, $A(x) - A(\bx)= 0$, and so
$$
\eqalign{-\int_x^\bx u(\eta)\; d\eta &=
-h(x, \beta(x)) + h(\bx, \beta(\bx)) - 
\int_{\beta(x)}^{\beta(\bx)} u(\eta) \; d\eta\cr
&= -h(x, \beta(x)) + h(x, \beta(\bx)) - 
\int_{\beta(x)}^{\beta(\bx)} u(\eta) \; d\eta\cr
&\qquad -h(x, \beta(\bx)) + h(\bx, \beta(\bx))\cr
& = \int_{\beta(x)}^{\beta(\bx)} (h_2(x,\eta) - u(\eta)) \; d\eta
+ \int_x^\bx h_1\bigl(x, \beta(\bx)\bigr) \; d\eta\cr
& >- \int_x^\bx u(\eta)\; d\eta\cr }$$
where in the last inequality we have used
the fact that $\beta$ is strictly increasing. \QED

\medskip
\noindent{\bf \S 2 Convex curves and envelope coordinates.}

In this section we review some standard material
from the theory of convex curves in the plane. For more
details the reader is referred to [E], [G1], and [G2].

\medskip
\noindent{\bf \S 2.1 The height function and the radius of curvature.}
If $\Gamma$ is a convex closed curve in the plane,
let $U(\G)$ denote the
unbounded component of $\reals^2-\G$. It will be convenient to
assume that the origin is {\it not} contained in $U(\G)$ and that
\mG\ is oriented in a counterclockwise direction.
Let $\vu_\theta = (\cos(\theta), \sin(\theta))$
and $\vn_\theta = (-\sin(\theta), \cos(\theta))$.
For each $\theta\in S^1 = \reals/(2 \pi \integers)$, $\LL_\theta$ is 
the supporting line to \mG\ that is parallel to
$\vn_\theta$ and near its point of intersection \mG\ is oriented in the
same direction as $\vn_\theta$. It will also be
convenient to assume that $\LL_0$ intersects \mG\ in an extremal point.
The {\de height function} of \mG, $p(\theta)$,
is the distance from the origin to $\LL_\theta$ (see Figure 2.1(a)). 
Equivalently, 
$p:S^1\ra \reals$ is given by $p(\theta) = G(\vu_\theta)$ where
$G(z) = \max\{z \cdot w : w\in \G\}$.

If $\LL_\t$ intersects \mG\ in an extremal point,
let $\va(\theta)$ be this intersection.
If $\LL_\t$ intersects \mG\  in a non-trivial segment, 
let $\va(\theta)$ be the extreme endpoint of the intersection 
in the counterclockwise direction. 
The arc length along \mG\ from $\va(0)$ to $\va(\theta)$ is denoted by 
$s(\theta)$. Note that  $s$ is  nondecreasing,  continuous
from the right and 
$$\va(\theta) - \va(0) = \int_0^\theta \vn_\eta\;  ds(\eta) $$ 
as a Riemann-Stieltjes integral. The points
of discontinuity of $s$ are exactly the countable set of  $\t$ for which
$\LL_\t$ intersects \mG\ in a non-trivial segment.

\longfig{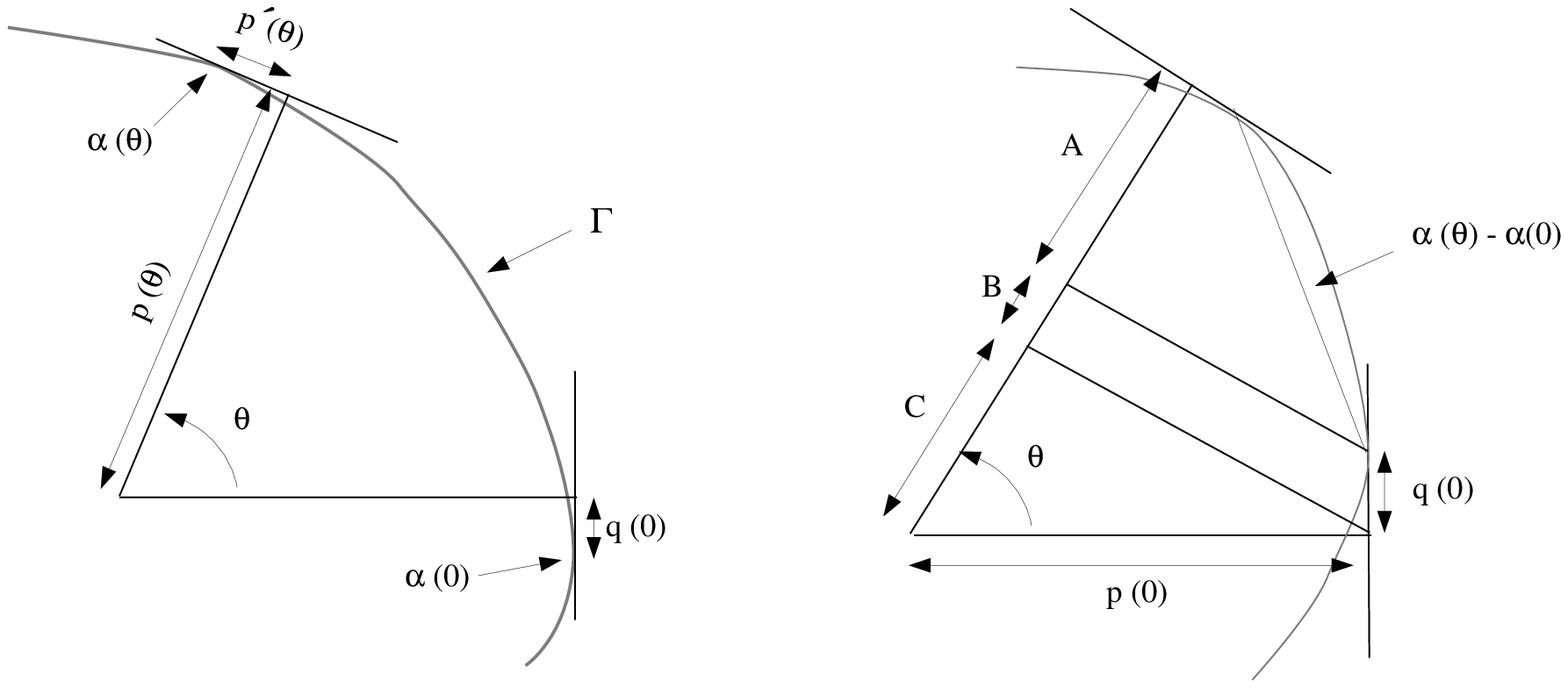,height=.40\hsize}{2.1}{(a) The height function; (b)
The geometry of equation (2.1)}

 A simple geometric argument (see Figure 2.1(b))  shows that
$$\eqalignno{p(\theta) & = 
-\vn_\t\cdot (\va(\theta) - \va(0))
+q(0) \sin(\t) + p(0) \cos(\t) \cr
& = \int_0^\theta \sin(\t - \eta) ds(\eta) 
+q(0) \sin(\t) + p(0) \cos(\t).&(2.1)\cr}$$
In Figure 2.1(b) the three terms in this expression
are labeled A, B and C, respectively. The quantity 
$q(0)$ is as pictured, explicitly, $q(0)$ is the signed distance
from $\alpha(0)$ to $p(0)\, \vu_0$   (the $q(0)$ shown in Figure
2.1(b) is positive). Note that the last two terms in (2.1) just reflect
the position of the origin.

It follows from (2.1) that $p$ is differentiable at all but
a countable number of points where $p^\prime$
has a simple jump discontinuity and 
$$p^\prime(\t) = \int_0^\theta \cos(\t - t) ds(t) 
- p(0) \sin(\t) + q(0) \cos(\t).\eqno(2.2)$$
At points where $p^\prime$ has a jump discontinuity
(those $\t$ for which $\LL_\t$ intersects \mG\ in a segment)
this formula may be interpreted in terms of right and
left hand limits.
In particular, $q(0) = p^\prime(0)$ and $p^\prime$ has a geometric
interpretation  as the signed distance
from   $\alpha(\t)$ to $p(\t)\, \vu_\t$ (see Figure 2.1(a)).
We shall be primarily interested in the case when \mG\ is
strictly convex which
corresponds to the condition that $s(\t)$ is continuous and 
$p$ is $C^1$.

If the arclength function $s(\t)$ is differentiable, 
$\rho(\t) := ds/d\t $ is the 
radius of curvature of \mG, \ie\ the radius of the osculating
circle to \mG\ at the point $\alpha(\t)$.
Using (2.1) and (2.2),  
$$p^{\prime\prime} + p = \rho \ge 0.\eqno(2.3)$$
If the $s(\theta)$ is not differentiable, this equation
may be interpreted in terms of distributions. For example,
if \mG\ is a polygon, $\rho$ is the sum of weighted Dirac
delta functions. In any case, since $s$ is  nondecreasing,
$\rho$ exists almost everywhere.
One interpretation of (2.3) is that the height function is
given by  a forced harmonic oscillator with forcing function
given by the radius of curvature (in which case the last two terms in
(2.1) reflect the initial conditions). Very informally, a closed
curve is a point that has been forced outward by the radius of
curvature. This observation will be crucial in connecting
\DB\ with impact oscillators in \S 7.

\medskip

{\bf Remark 2.1.1:} For later reference we give a characterization due 
to J. Green of the type
of periodic functions that can arise as 
height functions for a convex closed curve in the plane
([G1], [G2]).
The characterization involves a generalized convexity.
Let $S(x; A,B) = A \cos x + B \sin x$.
A function $g:\reals\ra\reals$ is said to be {\de sub-sine } if for
each $x_1 < x_2< x_1 + \pi$, $S(x; A_0, B_0) \ge g(x)$ for all
$x\in [x_1, x_2]$, where $A_0$ and $B_0$ are the unique values
for which $S(x_1;A_0, B_0) = g(x_1)$ and   
$S(x_2;A_0, B_0) = g(x_2)$. Green shows
that a $2 \pi$-periodic function $p$ is the height function of some convex closed curve if and only if $p$ is sub-sine. 
This happens if and only if $p$ can be written as
in (2.1) for some nondecreasing function $s$ with
$\int \exp( i \t) \;ds(\t) = 0$.  
Further, a twice differentiable $p$ is sub-sine if and only if 
$p^{\prime\prime} + p \ge 0$.

In a similar vein note that if $s :[0, 2 \pi)\ra \reals$  
is any bounded, nondecreasing
function with $\int \exp(i \t) \;ds(\t) = 0$
we can construct a \mG\ using (2.1) and (2.2).
The resulting curve will have 
parameterization $\va(\t) = R_\t (p(\t), p^\prime(\t -))$ where
$R_\t$ is rotation in the plane by the angle $\t$ and
$p^\prime(\t -)$ is the derivative at $\t$ from the
left. The curve
\mG\ will be the convex hull of $\va(\circle)$.
One can also use (2.1) to identify 
convex curves whose  perimeter has length one 
with regular Borel probability measures $ds$ on $S^1$
whose first Fourier coefficient vanishes. 
Atoms in the measures give flat spots (intervals in \mG);
strictly convex curves correspond to measures without atoms.

\medskip

\noindent{\bf \S 2.2 Envelope coordinates.}
The next step is to use \mG\ to give coordinates to the
open (topological) annulus $U(\G)$.  
For a point given in Euclidian coordinates  $z\in U(\G)$, 
let $\t(z)$ be the unique angle $\t$ with $z\in \LL_\t$
and $\va(\t)- z$ is parallel to 
$\vn_\t $. A second coordinate for $z$
is given by the distance from $z$ to $\va(\t)$ along $L_{\t(z)}$
and is denoted $\ell(z)$ (see Figure 2.2). Note that the area elements are
related by $dz = \ell\; d\ell\; d\t = d\gamma \; d\t$ where
$\gamma = \ell^2/2$. Thus the map
$z\mapsto (\t, \g)$ is area preserving. 
The coordinates $(\t, \g)$ 
will be called {\de envelope coordinates}. If \mG\ has corners, \ie\
its radius of curvature function $\rho(\t)$ vanishes on a nontrivial
interval, the map $z \mapsto (\t, \g)$ cannot be extended to a single
valued function on \mG.

\fig{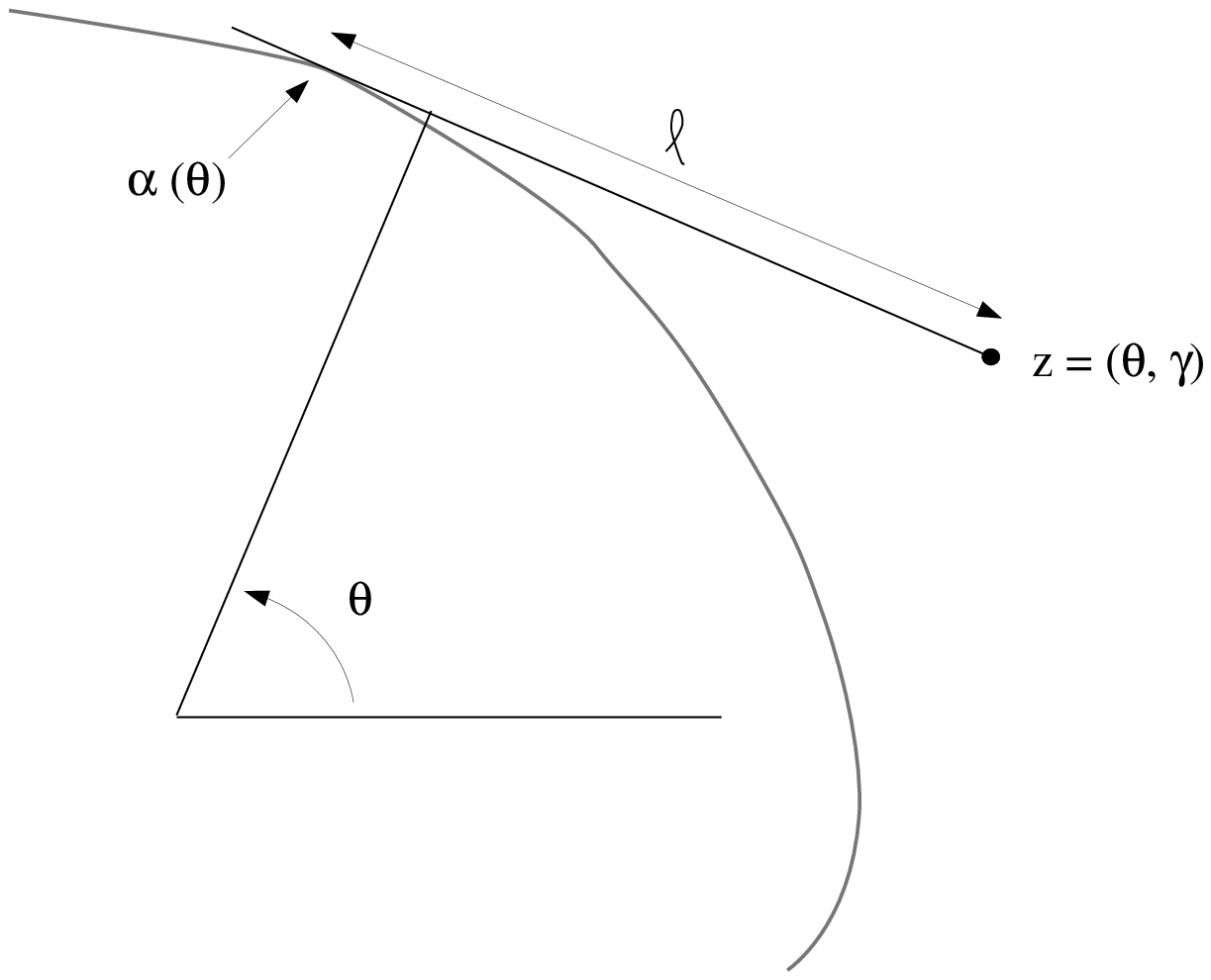,height=.30\hsize}{2.2}{Envelope coordinates}

\medskip

\noindent{\bf \S 3 Dual billiards yields a half-cylinder twist map}

A geometric description of the dual billiard map arising from
a convex curve \mG\ was given in the introduction.
The purpose of this section is to show
that under appropriate assumptions on
\mG\ the dual billiards map in
envelope  coordinates is a half-cylinder twist map.
It is geometrically clear that a dual billiards
map preserves Euclidian area. The corresponding map
in envelope coordinates is area-preserving because
 the change to envelope coordinates preserves area.  
The other main property that characterizes 
 twist maps is that the image of a vertical ray is a graph over 
an arc in the circle.
 A vertical ray above the
point $(\t_0, 0)$ in envelope
coordinates corresponds to the negative ray in
$\LL_{\t_0}$ beginning at the point of tangency with \mG. The image
of this ray under the \DB\ map  is the positive ray in
$\LL_{\t_0}$ that begins at the same point. Thus the
condition that makes $f$ a twist map is that this ray should
hit each other $\LL_\t$ ($\t_0 < \t < \t_0 + \pi$) in exactly one point. 
This follows from the convexity of \mG.

The dual billiards map on $U(\G)$ is denoted
$\Phi_\G$, or just $\Phi$, and
the corresponding map in envelope coordinates is denoted
$f_\G$ or $f$, and its lift to $\tC$ as $\tf$. 

\medskip

\noindent{\bf \S 3.1 The recursion relation from geometry.}
Let \mG\ be a convex curve in the plane. It has associated
functions $\alpha$, $p$, etc. as defined in \S 2.1. 
Given two angles $\t < \tp$, let $\delta = \tp -\t$ and 
$\phi = \a(\tp)- \a(\t)$. If $R = R(\t, \tp)$, 
$L = L(\t, \tp)$, $\omega_1$, and $\omega_2$
are as pictured in Figure 3.1, then the law of sines yields
$${\sin(\pi - \delta)\over |\phi|} = 
{\sin(\omega_1)\over L} =
{\sin(\omega_2) \over R}.\eqno(3.1)$$

\fig{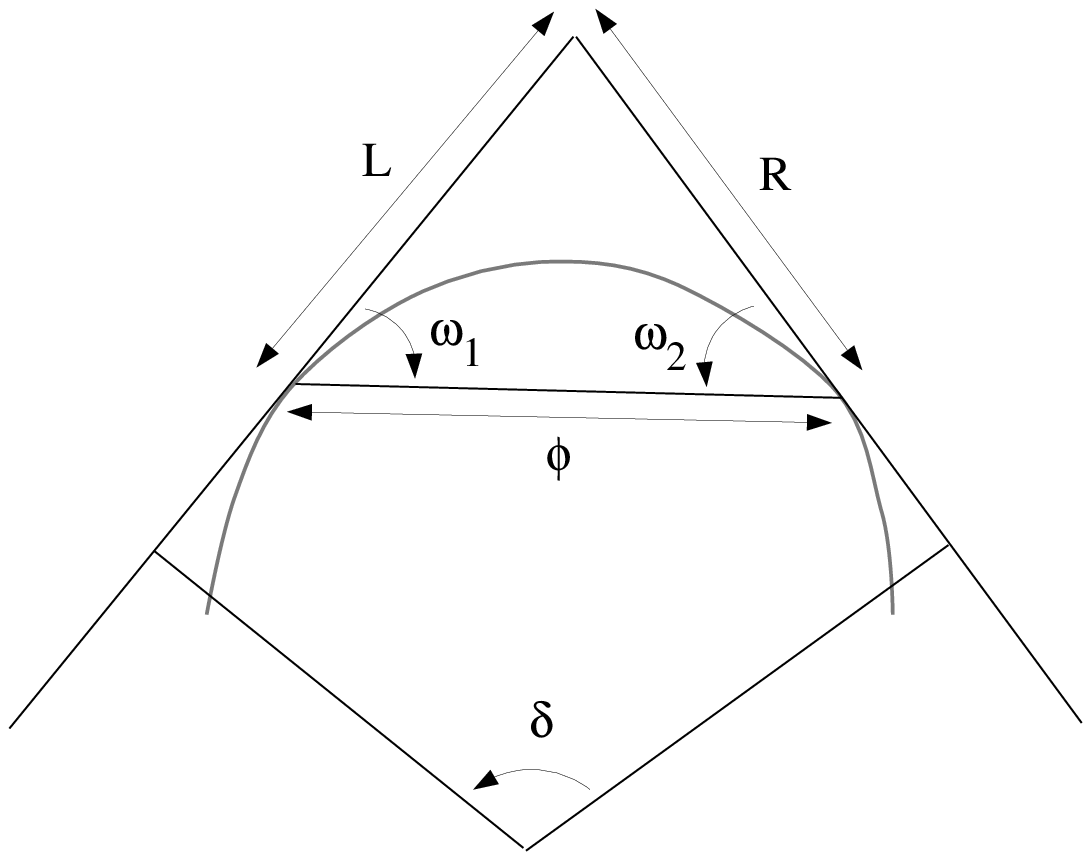,height=.30\hsize}{3.1}{Geometric basis of the recursion relation}

The wedge (or cross product) between two vectors in the
plane is ${\bf w}\wedge {\bf v} = w_1 v_2 - v_2 w_1 =
|w| |v| \sin(\omega)$ where $\omega$ is the angle from
${\bf w}$ to ${\bf v}$.
Using this notation and (3.1) we have
$$ L(\t, \tp) = {|\phi| \sin(\omega_1)\over \sin(\pi - (\tp - \t))}
= {\vn_\t \wedge \int_\t^\tp \vn_\eta\d s(\eta) \over \sin(\tp-\t)}
= {\int_\t^\tp \sin(\eta-\t)\d s(\eta) \over \sin(\tp-\t)},\eqno(3.2)$$
and a similar calculation yields
$$R(\t, \tp) = {\int_\t^\tp \sin(\tp - \eta)\d s(\eta) 
\over \sin(\tp-\t)}.\eqno(3.3)$$
Thus using the geometric definition of the dual billiard map
a triple of points $ (\t, \tp, \t^{\prime\prime})$ are the angles
of an orbit of dual billiards if and only if $0<\tp-\t < \pi$,
$0< \t^{\prime\prime} - \tp <\pi$, and
$$L(\t, \tp) = R(\tp, \t^{\prime\prime}) > 0. \eqno(3.4)$$
\medskip

\noindent{\bf \S 3.2 The generating function and the boundary maps.}
Next we define the functions $h$ and $b_i$ that are the
generating function and boundary maps for $f_\G$.
If we let
$$
h(x, \xp) = {1 \over 2} \int_x^{\xp} \bigl(L(x, y)\bigr)^2 \;dy
= {1 \over 2} \int_x^{\xp} \bigl(R(y, \xp)\bigr)^2\; dy\eqno(3.5)
$$
then $h_1 = - R^2/2$, $h_2 = L^2/2$,
and  the condition that $h_1(x^\prime, x^{\prime\prime}) +
h_2(x, \xp) = 0$ is equivalent to (3.4).

To define the boundary maps,
let $b_0(x) = \inf\{z : z > x \ \hbox{and}\ \rho(z) >0\}$. Perhaps
put more simply,  $b_0(x) = x$ if $x$ is not in an interval in which
$\rho$ identically vanishes. If $x$ is in such an interval,
$b_0(x)$ is the right endpoint of that interval. 
The function $b_0$ is clearly continuous from the right, nondecreasing,
satisfies $b_0(x + 2 \pi) = b_0(x) + 2 \pi$
and has rotation number equal to zero. Let $b_\infty(x) = x + \pi$.

Let us assume now that  \mG\ is a convex curve whose
arclength function $s(\t)$ is Lipschitz. In this case the
radius of curvature  $\rho$ exists almost everywhere, is
bounded and satisfies  $\rho \ge 0$. The derivatives
at points where $\rho$ exists are:
$$\eqalign{
R_1(x, \xp) &= -\rho(x) + \cot(\xp -x) R(x, \xp)\cr
R_2(x, \xp) &= {L(x, \xp) \over \sin(\xp-x)}\cr
L_1(x, \xp) &= -{R(x, \xp) \over \sin(\xp-x)}\cr
L_2(x, \xp) &= \rho(\xp) - \cot(\xp -x) L(x, \xp).\cr}\eqno(3.6)
$$
It is now a simple matter to check that $h$ and the $b_i$ satisfy the
necessary properties to define a  \HTM. 

\medskip

{\bf Theorem 3.1:} {\it If \mG\ is a convex curve whose
arclength function $s(\t)$ is Lipschitz, then its dual
billiard map $f_\Gamma$ in envelope coordinates
is a half-cylinder twist map with rotation set equal
to $[0, 1/2]$.  }
\medskip

{\bf Remark 3.2.1:} The function $h(\t, \tp)$ has the geometric
interpretation
as the area bounded by \mG\ and the supporting lines $\LL_\theta$ and
$\LL_{\tp}$ ([Ms4])  (see Figure 3.2). This area is the same as the
geometric interpretation of the generating function
given in \S 1.1 because the change from
Euclidian to envelope coordinates preserves area,
 The variational formulation also has a nice
geometric interpretation in \DB. For example, to find
a period three orbit of the \DB\ map, one finds the
circumscribed triangle of smallest area (see Figure 3.3, left).
Theorem 1.1(a)  says that we can find a monotone periodic
orbit of all rotation numbers. This  corresponds to finding
circumscribed polygons of each ``rotation type'' (see Figure 3.3,
right).
Further, there are the irrational analogs of these
sets for each irrational in the rotation set.

\fig{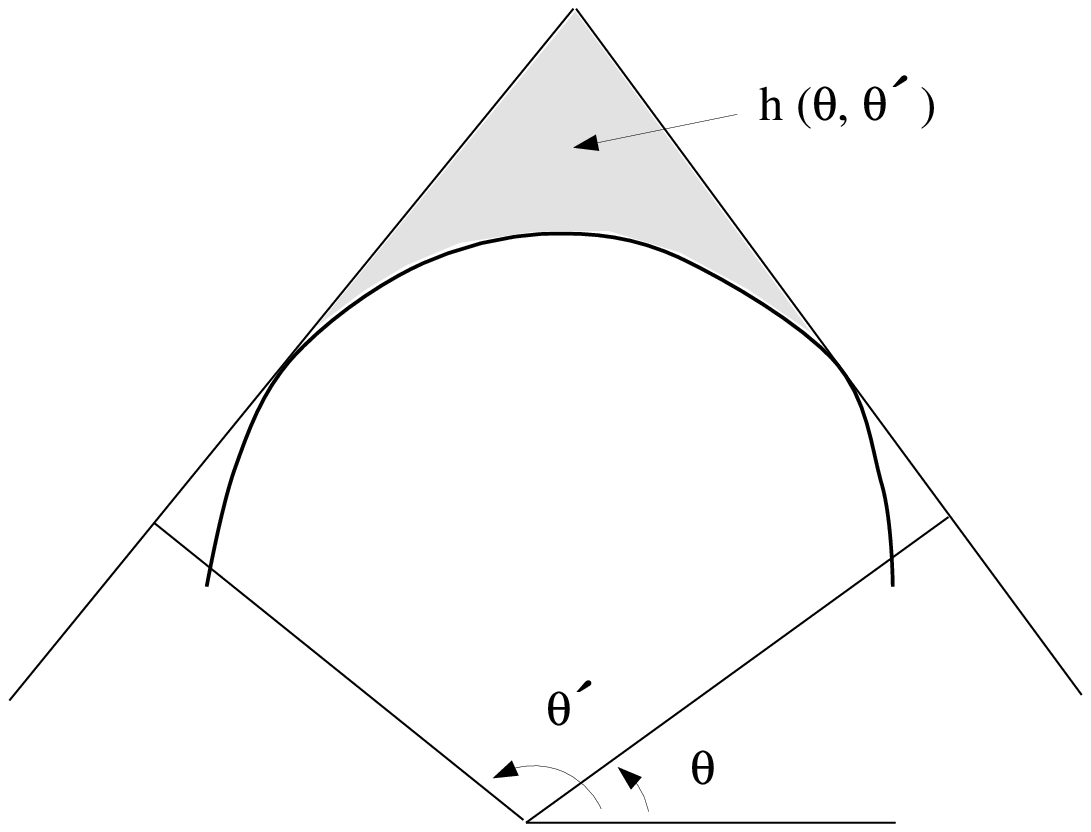,height=.30\hsize}{3.2}{The geometric interpretation of the action}

\longfig{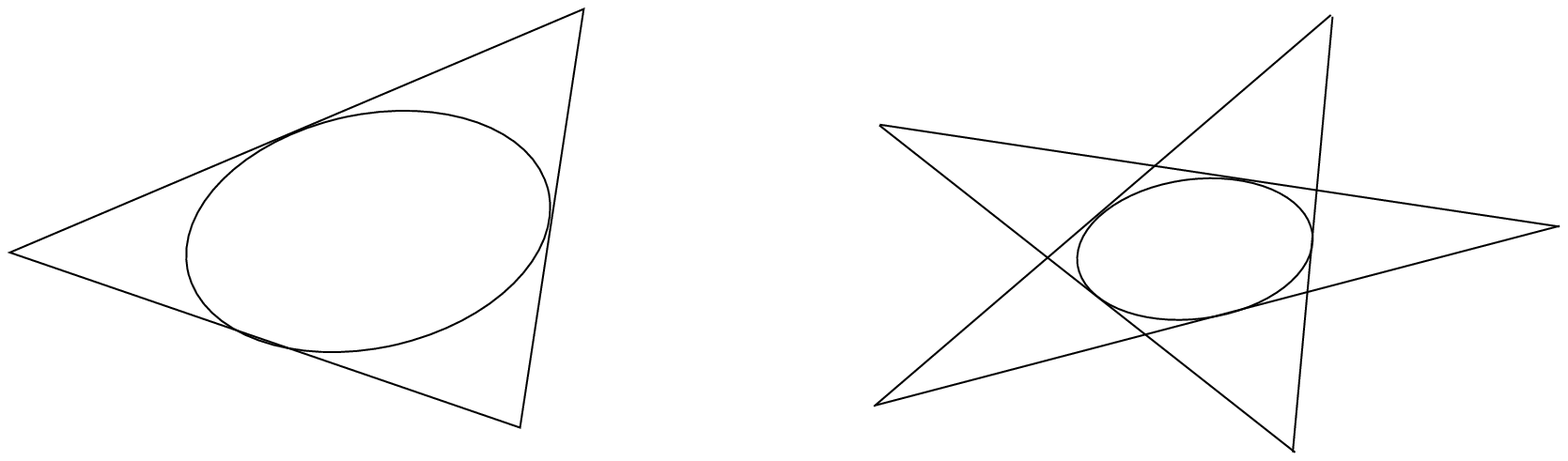,width=.60\hsize}{3.3}{Using variational methods to find monotone
periodic orbits with rotation number $1/3$ (left) and $2/5$ (right)}

{\bf Remark 3.2.2:} If the curve \mG\ is strictly convex, then its
arclength function $s(\t)$ is continuous, and the \HTM\
$f_\Phi$ on $C$ is also. By inspecting (3.6) we
see that this is sufficient to make
$h_{12}$  continuous, even though $h_{11}$, etc.
may not exist. If $s$ is $C^r$, then so is $f_\G$.

If $s$ is Lipschitz, $Dh$ will be also, and
so $f_\Phi$ is locally \Lip.
However, 
as $|x-\xp|\ra 0$, $h_{12}\ra 0$, and so (1.2) yields
that 
$\tf$ is never \Lip, no matter how smooth $s$ is.
 On the other hand, $\tf$ has bounded twist near 
the bottom boundary because 
$h_{11}$ and $h_{22}$ are bounded there. In fact,
as $|x-\xp|\ra 0$, $h_{11}\ra 0$, so the slope
of the image of a vertical tangent goes to zero
near the bottom edge.

Note that $\Phi_\G$ always extends continuously
to \mG\ by making it the identity there,  but
$f_\G$
will not have a continuous extension to
$\circle\times \{ 0 \}$ if \mG\ has
corners (nontrivial intervals where
$\rho$ vanishes). This is why discontinuous
boundary maps $b_0$ had to be allowed in
the definition of \HTM.
Note also that  as $\xp-x \ra \pi$ (\ie\ as 
$\g \ra \infty$), $h_{22}\ra\infty$ so the image
of vertical tangent is near vertical and the twist
of $f_\G$ goes to zero.

\medskip

\noindent{\bf \S 4 Invariant circles for dual billiards.}

As noted in the introduction, the question of the existence of invariant circles
is central in twist map theory. An invariant circle for a dual billiards map
is geometrically connected to the 
given convex curve via an equal area construction.
Recall that in this paper the phrase ``invariant circle'' always means
a homotopically nontrivial circle.
\medskip

\noindent{\bf \S 4.1 Existence of invariant circles.}
The first theorem  gives sufficient conditions on the 
curve $\G$ so that its \DB\ map has invariant
circles near to $\Gamma$ and near infinity.
The Hausdorff distance between two compact sets, $X$ and $Y$,
is denoted $d(X,Y)$.

\medskip
{\bf Theorem 4.1:} (R. \DD) {\it Let \mG\ be a convex curve in the plane
with a radius of curvature function $\rho(\t)$ that is $C^r$.

{\leftskip=30pt\rightskip=20pt\parindent=-18pt

(a) If $r>4$, then for each $N>0$ the set of $\Phi_\G$-invariant
circles $\Omega$ with  $d(\Omega, \G) > N$ has positive
\Leb\ measure.

(b) If $r>5$ and $\rho > 0$, then for each $\epsilon > 0$
the set of $\Phi_\G$-invariant circles $\Omega$ with 
 $d(\Omega, \G) < \epsilon$ has positive \Leb\ measure.

}} 
\medskip

The basic idea of the proof in [D] (see also [Ms2]) is to extend the
\HTM\ $f$ associated with $\Phi_\G$ to the boundaries
of $\Ca$ and then use a normal form result based on
Herman's version of the curve translation theorem ([H]).
In the compactification of the 
half-cylinder, $f$ does not preserve a finite measure near infinity.
However, it does satisfy the circle intersection property,
so the curve translation theorem can be used to obtain
invariant circles. The condition $\rho > 0$ in
(b) is needed to extend $f$ smoothly
to  $\reals\times \{0\}$ (\cf\ Theorem 4.3(a) below).

\medskip

\noindent{\bf \S 4.2 Invariant circles and area envelopes.}
The next result concerns 
the inverse problem for invariant
circles, \ie\ given a closed curve $\G_1$ can
you find a convex curve $\G_0$ so that  $\G_1$ is an
invariant circle for $\Phi_{\G_0}$.
The solution involves a geometric construction
called the  area envelope. 

Given a simple closed curve $\G_1$ in the plane 
and a number $a$ less than half the area enclosed by $\G_1$,
for each angle $\t$ we can find a unique oriented line
$\LL_\t$  
in the direction $u_\t$ so that the
area to the right of the line and inside $\G_1$ is equal to $a$.
Denote the  envelope of these lines as $AE(\G_1, a)$ (see Figure
4.1). If $AE(\G_1, a)$ is convex, 
it is an easy geometric exercise to show that
the locus of the midpoints of the
chords of the lines $\LL_\t$ in $\G_1$ is the area envelope.\footnote{%
$^\sharp$}{The author learned
this from R. Ticciati who ``discovered'' it during a High School project
in 1967}
The next result is a consequence of this geometric
fact. It appeared first in [D].

\fig{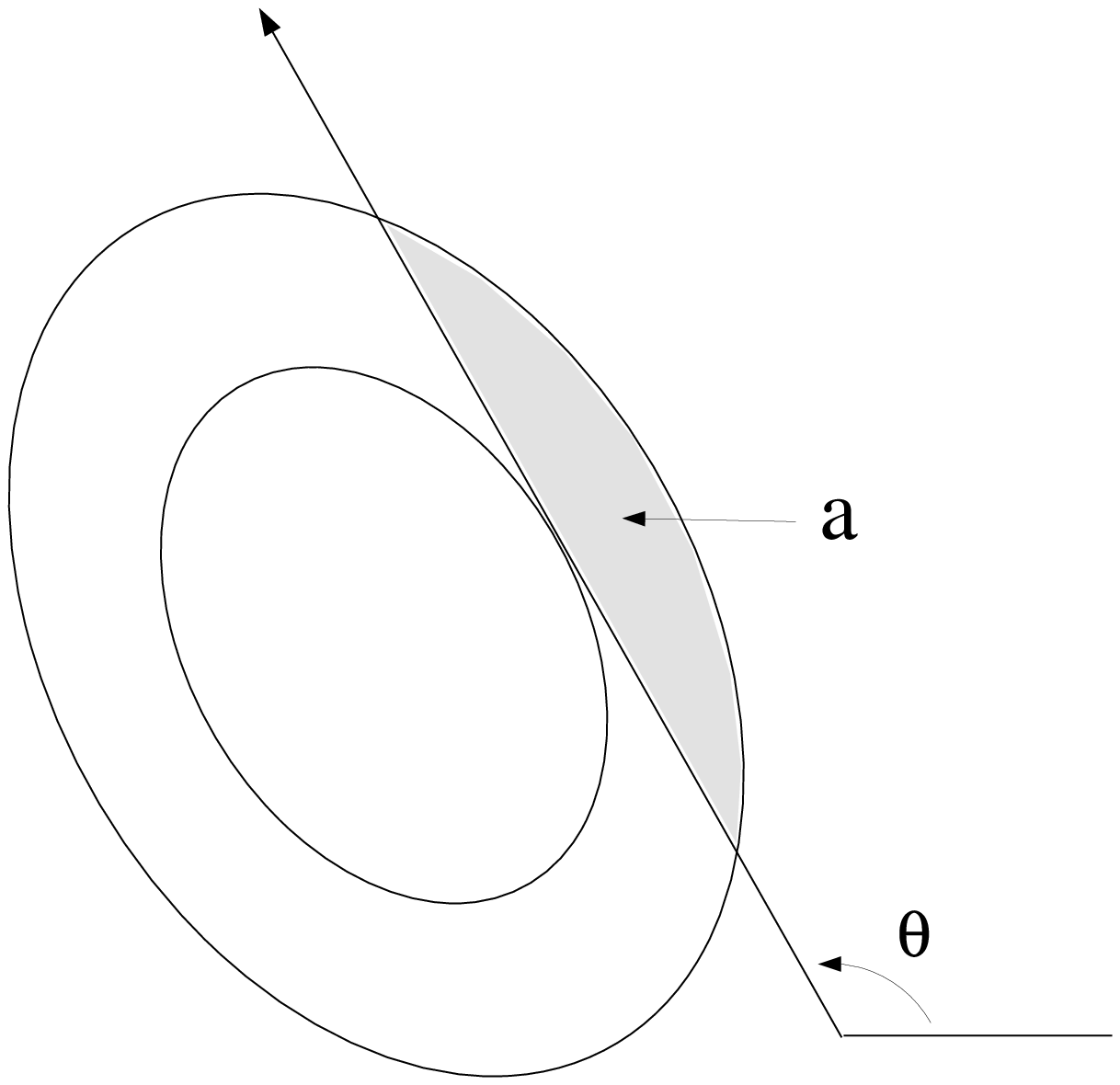,height=.30\hsize}{4.1}{An area envelope}

\medskip
{\bf Theorem 4.2:}  (R. \DD) {\it 
If $\Gamma$ is a convex curve with a \Lip\
arclength function,
then a closed curve $\Gamma_1\subset U(\G)$ is invariant
under the dual billiards map $\Phi_\G$ if and only
if $\Gamma = AE(\Gamma_1, a)$ for some $a$.
}
\medskip

Recall that the change to envelope coordinates is
area-preserving. Thus 
if one considers the \HTM\ corresponding to
$\Phi_\G$, the quantity $A(x)$ (defined in \S 1.4)
associated to  the invariant circle 
is exactly the appropriate area  $a$ for which
$\Gamma = AE(\Gamma_1, a)$.
Thus Theorem 4.2 is a consequence of the more general fact 
given in Proposition 1.3. The appendix contains a remark
on an analogous result for inner billiards.

Even if $\G_1$ is convex, 
its area envelope need not be (\eg\ an area envelope
inside an equilateral triangle). With the appropriate
interpretation of chord, the area envelope will
still be the locus of the midpoints.  Area envelopes
are discussed in more detail in [FT] and [GK].
 
\medskip
\noindent{\bf \S 4.3 The non-existence of invariant circles.}
The next result gives sufficient conditions
on \mG\ that imply the non-existence of invariant
circles near \mG. These results are the analogs for dual
billiards of results for billiards due
to Mather and Hubacher.  Conditions
on the curvature in billiards play the role of the conditions
on the radius of curvature in dual billiards.

The proof of part (a)  goes exactly like that
of Mather in [M1], and indeed, Theorem 1.1(c) is taken
from that paper. The proof of (b) is unlike that
of [Hr] in that we use estimates involving
the recursion relation of the twist map, rather than 
geometric arguments.
Part (a) is also contained in [GK]. That paper also
contains estimates  based on the geometry of \mG\
for the sizes of regions in $U(\G)$ 
that do not contain invariant circles.

\medskip
{\bf Theorem 4.3:} {\it Let $\Gamma$ be a convex curve in
the plane  with Lipschitz
arclength function $s(\t)$, and $J\in \circle$
is an open interval. If

{\leftskip=30pt\rightskip=20pt\parindent=-18pt

(a) the arclength function $s(\t)$ is $C^1$ on $J$, $\rho>0$
at the endpoints of $J$,  and for some
$\htheta\in J$, $\rho(\htheta) = 0$,

}

\noindent or

{\leftskip=30pt\rightskip=20pt\parindent=-18pt

(b) there is a $\htheta\in J$ so that $s$ is  $C^1$ on $J - \htheta$,
$\rho > 0$ on $J$, and $\rho$ has a jump discontinuity
at $\htheta$,

}

\noindent then there exists a neighborhood of $\G$ that contains no
homotopically nontrivial $\Phi_\G$-invariant circles.}
\medskip

{\bf Proof:} The proof will use the half-cylinder twist
map $f$ associated with $\G$ that has lift $\tf$.
Let $h$, $b_0$ and $\cD$ be the
generating function, lower boundary map and the domain of $h$, respectively,
and $\hat x$ is a lift of $\htheta$.

If $\tilde\Omega\subset\tC$  is the lift of
a homotopically nontrivial $f$-invariant circle, then by Theorem
1.1(b), $\tO$ is the graph
of a Lipschitz function $u:\reals\ra (0,\infty)$.
Using the derivative formulas (3.6) and Theorem 1.1(d), for almost every
$x_0$, 
$$2\rho(x_0) > \cot(x_1-x_0) R(x_0, x_1) + \cot(x_0-x\mo) L(x\mo, x_0)
\eqno(4.1)
$$
where $(x_i,\g_i) = f^i(x_0, u(x_0))$.

Assume first that the situation in (a) happens. If there were
 a $\Phi$-invariant circles arbitrarily close to 
to $\G$, then by Theorem 1.1(b),
we could find the lift of an invariant circle $\tO$ arbitrarily
close to $\reals \times \{ 0 \}$ in the Hausdorff topology, in which
case  $\max \{ |x_1 - x_0| , |x_0- x\mo| \} < \pi/2$ and
$\{x\mo,x_0, x_1\} \subset J$ for all $x_0$ sufficiently close
to $\hat x$.
In this case the right hand side of (4.1) 
is positive, continuous and bounded
away from zero while the left hand side is a continuous function
that vanishes at $x_0 = \hat x$. Thus the equation cannot hold
almost everywhere for $x_0$ in a neighborhood of  $\hat x$, a contradiction.

Now assume the situation given in (b). Let $\lim_{x\ra\hat x^\pm} \rho(x):=
\rho^\pm$ and assume without loss of generality that $\rho^+ > \rho^-$ and 
$|J| < \pi/2$. The principle estimate needed 
is the following: $\forall \epsilon>0$,
$\exists \delta_1 > 0, \forall N>0$, 
$\exists \delta_2 = \delta_2(\delta_1, N)$ so that $|x\mo-x_0|<\delta_1$,
$\delta_1/N < |x_0 - x_1 | < \delta_1$ and $0<\hat x - x_0 < \delta_2$
implies 
$$\max \{ |\bar L(x\mo, x_0) - \rho^-|, |\bar R(x_0, x_1) - \rho^+|,
|\rho(x_0) - \rho^-|\} <\epsilon$$
where $\bar L(x\mo, x_0) = \cot(x_0-x\mo) L(x\mo, x_0)$ and
$\bar R(x_0, x_1) = \cot(x_1-x_0) R(x_0, x_1)$.

To prove this first note that if $(x\mo, x_0, x_1)$ is a segment
of an  allowable configuration for $h$, 
$\{x\mo, x_0, x_1\} \in J$, and $\hat x > x_0$,
then
$$\eqalign{|\bar R(x\mo, x_0) - \rho^+|
&= \big| \cot(x_1-x_0) 
\int_{x_0}^{x_1} \bigl({\sin(x_1 - \eta)\over \sin(x_1-x_0)} \rho(\eta)
- {\tan(x_1 - x_0)\over (x_1 - x_0)}\rho^+ \bigr)\; d\eta \big|\cr
&\le {1\over (x_1 - x_0)} \int_{x_0}^{x_1} |\rho(\eta) - \rho^+|\; d\eta\cr
&\le {1\over (x_1 - x_0)} \bigl((\hat x - x_0)\sup\{|\rho(\eta) - \rho^+ | :
\eta\in [x_0, \hat x)\}\cr
&\qquad +  (x_1 - \hat x)\sup\{|\rho(\eta) - \rho^+ | :
\eta\in (\hat x, x_1]\}\bigr)\cr }$$

and
$$\eqalign{|\bar L(x\mo, x_0) - \rho^- |
&\le {1\over (x_0 - x\mo)} \int_{x\mo}^{x_0} |(\rho(\eta) - \rho^-)|\; d\eta\cr
&\le \sup\{ |\rho(\eta) - \rho^-| : \eta\in [x\mo, x_0]\}.\cr}$$

Now given $\epsilon > 0$, since $\rho$ is continuous and bounded
on $J - \hat x$ there is an $\delta_1>0$ so that $0 < x_1 - \hat x < 
\delta_1$ implies
$$\sup\{|\rho(\eta) - \rho^+ | : \eta\in (\hat x, x_1]\}  < \epsilon/2$$
and in addition, $0< \hat x - x\mo < 2 \delta_1$ implies
$$\sup\{ |\rho(\eta) - \rho^-| : \eta\in (x\mo, x_0)\} < \epsilon.$$

Now given $N$, pick $\delta_2 < \delta_1$ so that
$${N \delta_2 \over \delta_1} \sup\{|\rho(\eta) - \rho^+ | :
\eta\in [x_0, \hat x)\} < \epsilon /2.$$
With these choices,  note that 
$\hat x - x\mo = x_0 - x\mo + \hat x - x_0 < \delta_1 +
\delta_2 < 2 \delta_1$, and thus principle estimate follows.

Continuing the proof of the lemma, pick $\epsilon > 0$ so that
$\rho^+ - \rho^- > 4 \epsilon$ and let $\delta_1$ be as 
in the principle estimate. If there are $f$-invariant circles arbitrarily
near $\circle\times \{0 \}$, then 
there is a lift of an invariant circle for which 
$0 < x_i - x_{i-1} < \delta_1$ for all $x_i\in J$ where
$x_i = \pi_1(\tf^n(z))$ with $z$ an element of the  lift of the  invariant
circle.  Fix one such circle
$\Omega$. Since it contains no fixed points by Theorem 1.1(b)
there is an $N$ with $\delta_1/N < 
 x_i - x_{i-1} < \delta_1$ for all pairs on $\Omega$. Using this
$N$, find $\delta_2$ as in the principle estimate. Then the
principle  estimate along with (4.1) implies that for
almost all $x_0$, in particular for some $x_0$ with
$0< \hat x - x_0 < \delta_2$,
$0 > \bar L(x\mo, x_0) + \bar R(x_0, x_1) + 2 \rho( x_0) >
\rho^+ - \rho^-  - 4 \epsilon > 0$, a contradiction.\QED
\medskip

{\bf Remarks:}

{\bf 4.3.1:} Neither the condition of part
(a) nor (b) excludes the existence of
invariant circles  outside neighborhoods of \mG. Indeed,
the area envelope construction of \S 4.2 allows one
to construct an invariant circle for \DB\ on a convex curve whose
radius of curvature function has any prescribed {\it local} behavior.
In the case of billiards, the analog of part (a) (a point
where the curvature vanishes) excludes all invariant
circles, while the analog of part (b) does not.

{\bf 4.3.2:} The Birkhoff-Mather stability theorem
relates the non-existence of invariant circles
with the existence of orbits with certain limit behavior.
The theorem states that if $\Omega_1$ and $\Omega_2$ are invariant
circles for a twist map, then there are no other invariant
circles in the annulus bounded by $\Omega_1$ and $\Omega_2$
  if and only if  there is a point
$z$ whose  $\alpha$- and $\omega$-limit sets are contained in 
$\Omega_1$ and $\Omega_2$, respectively. Mather's result is proven in
[M4] for a specific class of twist maps, but the arguments
are quite robust and almost certainly apply to
\HTM s. 

Thus, for example,  when there are no
 $\Phi_\G$-invariant circles near $\G$ there is an
orbit that converges to \mG\ in forward time. This orbit is
distinguished from the crash orbits of the next section by the
fact that it converges to all of \mG, not a single point on \mG. There will also be an orbit that converges to $\G$ under
backward iteration, so it may be viewed as ``emerging'' from
\mG.

The question of stability  depends on  the behavior at infinity.
 In the Birkhoff-Mather theorem, one of the invariant circles can 
be located at infinity. Thus if there are no invariant circles near infinity,
there is an orbit that escapes, \ie\ goes
to infinity under forward iteration.
 As noted in the introduction
one of the most important outstanding problems in the theory
of \DB\  is the existence of a 
convex curve whose \DB\ map  has such an orbit.

{\bf 4.3.3:}  Remark 1.2.1 contains the main geometric idea underlying
the proof of Theorem 4.3. The conditions given in (a) or (b) ensure that
a vertical tangent tumbles past the vertical in two iterates. 
To obtain a stronger result on the non-existence of invariant circles, 
one would need to use more iterates. However, for 
the \HTM s coming from \DB, the twist goes to zero at infinity.
It thus requires more and more iterates to tumble a tangent as
one approaches  infinity. This partially explains the difficulty in 
obtaining a \mG\ that yields no invariant circles near infinity. 
As noted in Remark 4.3.1, no local condition on
\mG\ can eliminate these invariant circles.

\hfill\eject
\medskip
\noindent{\bf \S 5 Crash orbits}

A {\de crash orbit} for \DB\ is an orbit that converges under
forward iteration to a point on the convex curve  $\Gamma$. 
The first theorem states that there exist
convex curves with differentiable radius of curvature functions $\rho$
for which there are crash orbits. However, the second theorem says
that if  $\rho$ is positive and has bounded derivative, 
then there do not exist crash orbits. These are analogs
of billiards results due to Halpern ([Ha]).
The  method of proof here is  different, although the analog
of Hapern's construction would also suffice. 

We give a fairly explicit $\rho$ and then use the criterion
of \ANG\ from Theorem 1.2 to get the existence of
the crash solutions. 
The simplest version of the example
has a radius of curvature function near zero
that looks roughly like
$1 + (1/2) x^2 \sin (1/ x^2)$.  Note that this is a standard
example of a function whose derivative exist everywhere, but
the derivative is not bounded.

The existence of the crash orbits
implies that the \DB\ map has no invariant circles
near \mG.  However, as in Remark 4.3.1, one can construct examples
with crash orbits and invariant circles far away using the
area envelope construction.

\medskip
\noindent{\bf \S 5.1 Examples with crash orbits.}
\medskip
{\bf Theorem 5.1:} {\it There exists a  convex curve 
$\Gamma$ with a differentiable radius curvature function 
such that the corresponding \DB\ map has crash orbits.
}
\medskip

{\bf Proof:}
We use the \ANG\ criterion of Theorem 1.2  with a modified form of the 
recursion relation. 
Define $\overbar \Delta(x\mo, x_0, x_1) = R(x_0, x_1) - L(x\mo, x_0)$
where $R$ and $L$ are as in (3.2) and (3.3).
It is clear that solutions, subsolutions, etc. of this recursion
relation will be the same as those given by the $\Delta$ defined
in (1.5) with $h$ given by (3.5).

We will construct a subsolution $\underline c$
with $c_n \nearrow 0$, and a supersolution $\underline w$ with
$w_n \nearrow w < \infty$, so that $\underline c \le \underline w$.
By Theorem 1.2, there exists a solution $\ux$ with
$\underline c \le \ux \le \underline w$, and since any solution is
of necessity monotone increasing, $x_n$ converges as $n\ra\infty$.
Now  by (1.1), there exists
some $(x_0, \g_0)$ with $x_i = \pi_1(\tf^i(x_0, \g_0))$.
Then using (3.2) or (3.3), $\g_n\ra 0$, and so $\tf^i(x_0, \g_0)$
converges to a point on $\reals\times \{ 0\}$, thus the corresponding orbit
in $U(\G)$ converges to a point of $\Gamma$ under the \DB\ map.

Given a sequence $a_n$ (to be specified later)
with $n\in \naturals$ such that $a_n \nearrow 0$  
and $a_n\in [-.1, 0)$,
let $\delta_n = (a_{n+1}-a_n)/2$ and
$b_n = a_n + \delta_n$. Pick a $0 < c < .001$, $k>1$,
and for $n = 1, 2, \;\dots, $ define $\rho:\reals\ra\reals$ by
$$ \rho(x) = 
\cases{ 1- c a_n^k& for $x\in [b_{n-1}, a_n)$\cr
 1 + c a_n^k & for $ x\in [a_n, b_n)$\cr
1   & for $x\in [0, 2 \pi - b_0)$\cr} 
$$
and extend $\rho$ so that it is $2 \pi$-periodic.

Now let
$$r(x) = {\cos(x)- \cos(2x)\over \sin(2 x)}$$
and 
$$s(x)  = {1-  \cos(x)\over \sin(2 x)}.$$
Note that  at points where the denominator is
zero, $r$ and $s$ can be continued to real analytic functions
with $s(x) + r(x) = \tan(x)$.

If $r_n = r(\delta_n)$ and $s_n = s(\delta_n)$,
 a simple calculation yields that 
$$\eqalign{R_n &:= 
R(a_n, a_{n+1}) = (1 + c a_n^k) r_n + (1 - c a_{n+1}^k) s_n\cr
L_{n-1} &:= 
L(a_{n-1}, a_n) = (1 + c a_{n-1}^k) s_{n-1} + (1 - c a_n^k) r_{n-1}.\cr}
$$

Our goal is to show that 
$R_n - L_{n-1} \ge 0 $ for sufficiently large $n$. 
Assuming $\delta_{n-1} - \delta_n > 0$, this is equivalent to 
$$
\eqalign{ 
&c\; ( {a_n^k r_n - a_{n+1}^k s_n - a_{n-1}^k s_{n-1} + a_n^k r_{n-1}\over
\delta_{n-1} - \delta_n})\cr
&> {s_{n-1} + r_{n-1} - s_n - r_n\over
\delta_{n-1} - \delta_n}\cr
&= {\tan(\delta_{n-1}) - \tan(\delta_{n})\over \delta_{n-1} - \delta_n}\cr
}\eqno(5.1)
$$
which as $n\ra\infty$ goes to $\tan^\prime(0) = 1$.

Now specify $a_n = -n^b$ for some $b<0$. In this case,
$$ \eqalign{-b n^{b-1} &> \delta_n > -b (n+1)^{b-1}\cr
 b (b-1) (n+1)^{b-2} &< \delta_{n-1} - \delta_n < b (b-1) (n-1)^{b-2}.\cr}$$

Using Taylor's theorem, for sufficiently small positive
$\delta$, $r(\delta) > 3 \delta/4$ and $-s(\delta) > - \delta/2$.
Thus the left hand side of (5.1) is larger than
$$c( {{3\over 4} n^{kb}((n+1)^{b-1} + n^{b-1})-
{1\over 2} (n+1)^{kb} n^{b-1} - {1\over 2}(n-1)^{(k+1)b -1} \over
(1-b) (n-1)^{b-2} }),$$
which goes to infinity if 
 $k b + 1 > 0$. Thus for these choices of $k$ and $b$,
 $R_n > L_{n-1}$ for $n \ge N$ for some
sufficiently large $N$.

Now pick $\hat \rho$ that is $C^\infty$ except at $0$ and is close
enough to $\rho$ to ensure that $\underline a$ is still a subsolution
(for $n\ge N$)
for the recursion relation defined using $\hat \rho$.

To construct the subsolution $\underline c$, let
$\g_N = h_2(a_{N-1}, a_{N})$ and for 
$i\le 0$, $c_{N+i} = \tf^i(a_N, \g_N)$ and for 
$i\ge 0$, $c_{N+i} = a_{N+i}$. By construction,
$\Delta(c_{i-1}, c_i, c_{i+1}) = 0$ if $i \le N-1$ and
$\Delta(c_{i-1}, c_i, c_{i+1}) \ge 0$ if $i \ge N$, and thus
$\underline c$ is a subsolution.

To construct a supersolution $\underline w$ with $\underline w
> \underline c$ begin by letting $v_i = \pi_1(\tf^i(a_N, \g_N))$ where
$(a_N, \g_N)$ is defined in the previous paragraph. Now if
the sequence $v_i$ is bounded above it is a crash orbit and we are
done, otherwise let $i_0 = \inf\{ i : v_i > 0\}$ and $w_{N-j} =
v_{i_0 - j}$ for $j\ge 0$. Define $w_{N+j}$ for $j\ge 0$ inductively as
follows. Assuming $w_{n-1}$ and
$w_n$ have been obtained, pick $w_n < w_{n+1} 
 < w_N + \sum_{i=0}^{n+1-N}   {1 \over 2^{i + 3}}$
so that $0< \hat R(w_n, w_{n+1}) < \hat L(w_{n-1}, w_n)$, where
$\hat L$ and $\hat R$ are defined by (3.2) and (3.3) using
$\hat \rho$ for the radius of curvature. This is possible
because $\hat R(x, \xp)\ra 0$ as $|\xp - x| \ra 0$.

Now note that $\Delta(w_{k-1}, w_k, w_{k+1})$ is less than zero
for $k\ge N$ and equal to zero for $k<N$, and thus
$\underline w$ is a subsolution. By construction,
$\underline c \;\le\; \underline w$ and $\lim w_n \le w_N +  1/4$.\QED

\noindent{\bf \S 5.2 Nonexistence of crash solutions.}

\medskip
{\bf Theorem 5.2:} {\it If $\Gamma$ is a convex curve in the plane
with a radius of curvature function $\rho > 0$ that is differentiable
and the derivative is bounded, then the \DB\ map has no crash orbits.
}
\medskip
{\bf Proof:}
Assume to the contrary that the \DB\ map has a crash orbit.
Lift the corresponding \HTM\ to the universal cover
and say the crash orbit is $(x_n, \gamma_n)$, where we assume
without loss of generality that $x_n \nearrow 0$.

Using the recursion relation given by
(3.4) and the mean value
theorem for integrals we can find $\eta_{n-1}\in (x_{n-1}, x_n)$ and
$\tau_n\in (x_n, x_{n+1})$ with
$$\rho(\eta_{n-1}) \tan(\delta_{n-1}/2)
= \rho(\tau_n) \tan(\delta_n/ 2)$$
where $\delta_n = x_{n+1} - x_n$,
and thus
$$
\bigl( \rho(\eta_{n-1}) - \rho(\tau_n)\bigr) \tan(\delta_{n-1}/2)
= \rho(\tau_n)\bigl( \tan(\delta_n/2) - \tan(\delta_{n-1}/2).
$$

Now by construction, $\tau_n - \eta_{n-1} < \delta_n + \delta_{n-1}$
and if $\delta>0$ is small, $\delta < \tan(\delta) < 2\delta$. 
Using the mean value theorem for derivatives, there is
a $\sigma_n\in [\eta_{n-1}, \tau_n]$ with

$$\eqalign{ \big|\rho^\prime(\sigma_n)\big| &= 
\big|{\rho(\tau_n) - \rho(\eta_{n-1})\over \tau_n - \eta_{n-1}}\big|\cr
&= \big|{\rho(\tau_n) \bigl(\tan(\delta_n/2) - \tan(\delta_{n-1}/2)\bigr)
\over (\tau_n - \eta_{n-1}) \tan(\delta_{n-1}/2)}\big|\cr
&\ge \big|{\rho(\tau_n)\bigl(\delta_n/2 - \delta_{n-1}\bigr)\over
(\delta_n + \delta_{n-1})\; \delta_{n-1} }\big|\cr
&= \big|{\rho(\tau_n)\over 2}\bigl({{\delta_n\over\delta_{n-1}} -2 
\over \delta_{n-1} + \delta_n}\bigr)\big|.\cr}\eqno(5.2)
$$

Now since $\rho$ is continuous, as $n\ra\infty$,
$\rho(\tau_n)\ra \rho(0) > 0$;
since $\sum \delta_n$ converges, $\delta_n\ra 0$ and 
$\limsup |\delta_n/\delta_{n-1}|
\le 1$. Thus the last expression in (5.2) diverges, contradicting the 
boundedness of $\rho^\prime$.\QED
\medskip

\noindent{\bf \S 6 The normalized action of the derivative}

In this section we make some geometric observations about the 
action of the derivative of the \DB\ map on tangent directions.
Since our aims are primarily descriptive, we will adopt a more 
informal tone than in previous sections.

\medskip
\noindent{\bf \S 6.1 Computation of  the action of the derivative.}
 Let $\beta(t)= (\beta_1(t), \beta_2(t))$ be a
parameterization of \mG\ {\it by arc length},
and let us 
assume initially that $\beta$ is twice differentiable. 

Fix a point $z\in u(\G)$ and say its envelope
coordinates are $(\t, \g)$. Let $t$ and $l$ 
be such that $\beta(t) = \alpha(\t)$ and $\ell= \sqrt{2 \gamma}$. 
Thus $\beta(t)$ is
the  point on \mG\ used
for the \DB\ shot, and $\ell$ the distance from $z$ to
$\beta(t)$. 
From inspection of Figure 6.1 we have
$$\eqalignno{
2 \beta(t) &= z + \Phi(z) &(6.1)\cr
z &= \beta(t) - \ell \beta^\prime(t).&(6.2)\cr}$$

\fig{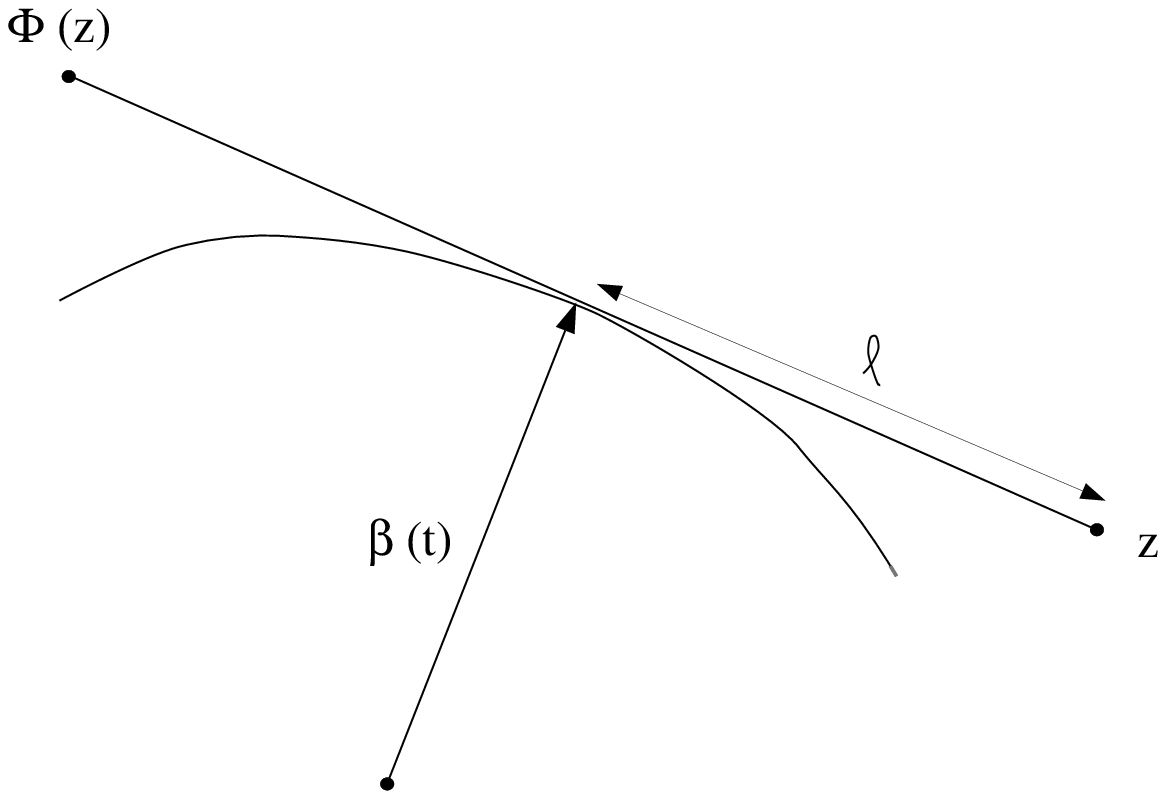,height=.30\hsize}{6.1}{Some geometry of \DB}

Adopting the classical (and often confusing) notation of
using the same symbol for a coordinate and a coordinate change,
  $z = (x, y)  = (x(t, \ell), y(t,\ell))$, and
subscripts to denote differentiation, 
(6.2) yields
$$ M := \pmatrix{x_t & x_l\cr y_t & y_l\cr } =
\left(\matrix{\beta_1^\prime - \ell \beta_1^{\prime\prime} &\beta_1^\prime\cr
\beta_2^\prime - \ell \beta_2^{\prime\prime} &\beta_2^\prime\cr}\right).
$$
Note that $det(M) = \ell(\beta^\prime \wedge\beta^{\prime\prime}) =
\ell \kappa(t)$, where $\kappa(t)$ is the curvature at the point
$\beta(t)$ (the curvature is the reciprocal of the radius of
 curvature $\rho$).
Thus inverting $M$ we get
$$
\left(\matrix{t_x& t_y\cr \ell_x&\ell_y\cr}\right) =
{1\over\kappa \ell}\left(\matrix{\beta_2^\prime&-\beta_1^\prime\cr
-\beta_2 + \ell\beta_2^{\prime\prime} & \beta_1 - \ell\beta_1^{\prime\prime}\cr
}\right).\eqno(6.4)
$$

Differentiating (6.1) yields
$$\eqalign{
D\Phi + Id &= 2 \pmatrix{\beta_1^\prime&0\cr\beta_2^\prime& 0\cr}
\pmatrix{t_x& t_y\cr \ell_x&\ell_y\cr}\cr
&= {2\over\kappa \ell}\pmatrix{
\beta_1^\prime\beta_2^\prime&-(\beta_1^\prime)^2\cr
(\beta_2^\prime)^2&-\beta_1^\prime\beta_2^\prime,\cr}
}
$$
and thus acting on a vector $\vu$,
$$
(D\Phi + Id)[\vu] = {2\over\kappa\ell}(\vu\wedge\beta^\prime)\;\beta^\prime.
\eqno(6.5)
$$

We now need some elementary geometry of the wedge product. Let
$\vw$ be a unit vector, $\vu$ and $\vv$ be vectors,
and $a,b,c$ be positive scalars so that 
$b \vw = a \vv +c \vu$.
If $A$ represents the area of the triangle formed by the vectors
then,
$$2 A = c\vu\wedge b\vw = b\vw\wedge a\vv = c\vu\wedge a\vv$$
and so
$$2 A = {b^2 (\vw\wedge\vv ) (\vw \wedge\vu)\over \vv\wedge\vu}.\eqno(6.6)$$
Now wedging (6.5) by $D\Phi[u]$ on the left
yields
$$
D\Phi[\vu]\wedge\vu = {2\over\kappa\ell}(\vu\wedge\beta^\prime)(D\Phi[\vu]
\wedge\beta^\prime)\eqno(6.7)
$$

%\fig{fig6_2.ps,height=.30\hsize}{6.2}{The triangle for equation (6.6)}

Thus if $B$ is the area of the shaded region in
in Figure 6.2, \ie\
$B$ is the area of the triangle bounded by the line connecting $z$ 
to $\Phi(z)$, the line in the direction  $\vu$ passing through the point $z$,
and the line in the direction of $D\Phi[\vu]$, then using (6.6) and (6.7),
$$
B = {(2\ell)^2\over 2}{(\beta^\prime\wedge D\Phi[\vu])(\beta^\prime\wedge\vu)
\over D\Phi[u]\wedge\vu}
= \kappa\ell^3.\eqno(6.8)
$$

Equation (6.8) contains the same information as Proposition 2.2
in [GK]. 

\fig{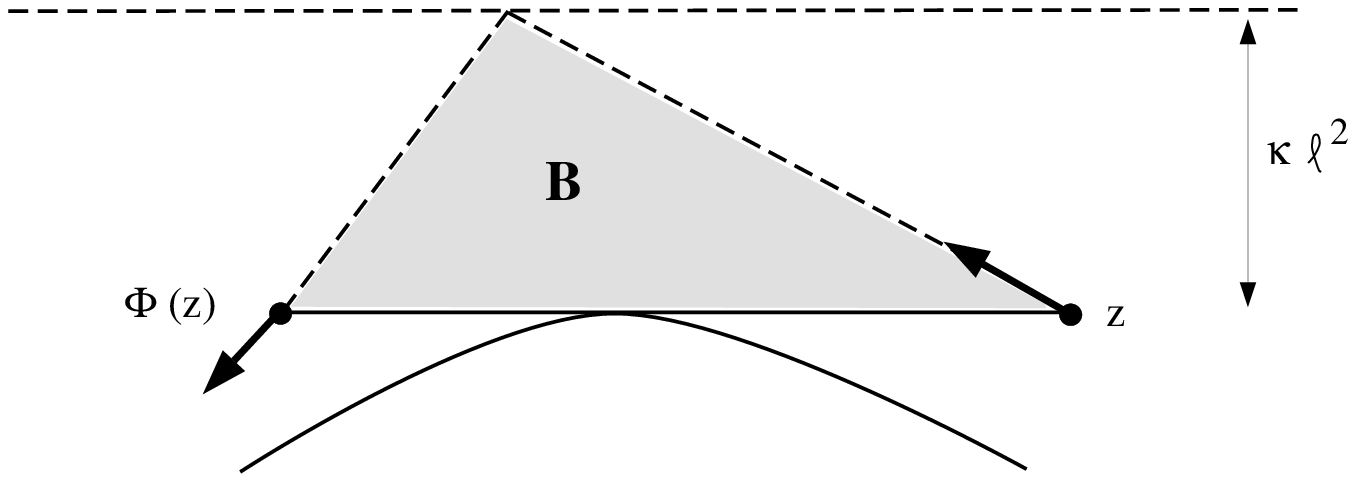,height=.25\hsize}{6.2}{Geometry of the derivative}

\medskip
\noindent{\bf \S 6.2 Geometric computation of the derivative.}
To geometrically compute
the action of the derivative of the dual billiards
map on a tangent direction, 
first fix a point $z$ and its image $\Phi(z)$ as in Figure 6.2. Now
draw a line parallel to that from $z$ to $\Phi(z)$ that is a distance
$\kappa \ell^2$ from that line, where $\kappa$ is the curvature at
the point of tangency. Call this new line the
{\de bounce line}.
Given a unit tangent vector
$\vv$ based at $z$, draw a line in the direction of $\vv$ until
it hits the bounce line. The line from this point of
intersection to the
point $\Phi(z)$ will give the direction of the vector $D\Phi[\vv]$.
This is because the shaded region in Figure 6.2 has area equal to $\kappa\ell^3$ as required by (6.8).

If \mG\ has a point with $\rho(\t) = 0$, by taking limits one 
sees that the action on tangent directions is a  bounce off infinity.
Thus when the dual billiards shot uses a corner of
\mG, the action of $D\Phi$ is rigid rotation by $\pi$.

{\bf Remark 6.2.1:} One motivation for seeking a geometric understanding 
the action of the derivative is to assist in the  finding
of strictly invariant
cone fields and thus (perhaps) obtain positive Lyapunov
exponents as in [W]. This theory does not seem to be well understood in
the case of infinite area as is the case in dual billiards.
One could also search for these cone fields in bounded
regions, \eg\ between a given curve and a second curve inside 
constructed using the  area envelope construction. The existence of positive 
exponents is  
particularly interesting in light of the physical 
interpretation of \DB\ given in \S 7. 

\medskip

\noindent{\bf \S 6.3 Comparison with a hyperbola.}
Having understood geometrically the action of the derivative we
are now in a position to interpret the tumbling tangent
criterion of Remark 1.2.1.  Again, pick a point $z$. Figure 6.3(a)
shows the action of the second iterate of a \DB\ map
on a tangent direction based at $\Phi\inv(z)$. There are three supporting lines
involved in the iteration, call them $\LL\mo, \LL_0$, and
$\LL_1$. The half-ray in $\LL\mo$ that begins at \mG\ and
points downward corresponds to a vertical line in
envelope coordinates. Thus the initial vector chosen
corresponds to the initial vertical tangent in Figure 1.2.
Its first iterate  under $D\Phi$ yields a parallel vector pointed
in the opposite direction. For the next
iterate we use the procedure of the previous subsection. Draw
a bounce line at height $\kappa \ell^2$ and compute the
next image. The condition that this vector has tumbled beyond 
the vertical  is that the image is on the 
 side of  $\LL_1$ nearest \mG. In the figure the initial vector has
tumbled by the third iterate. 

\longfig{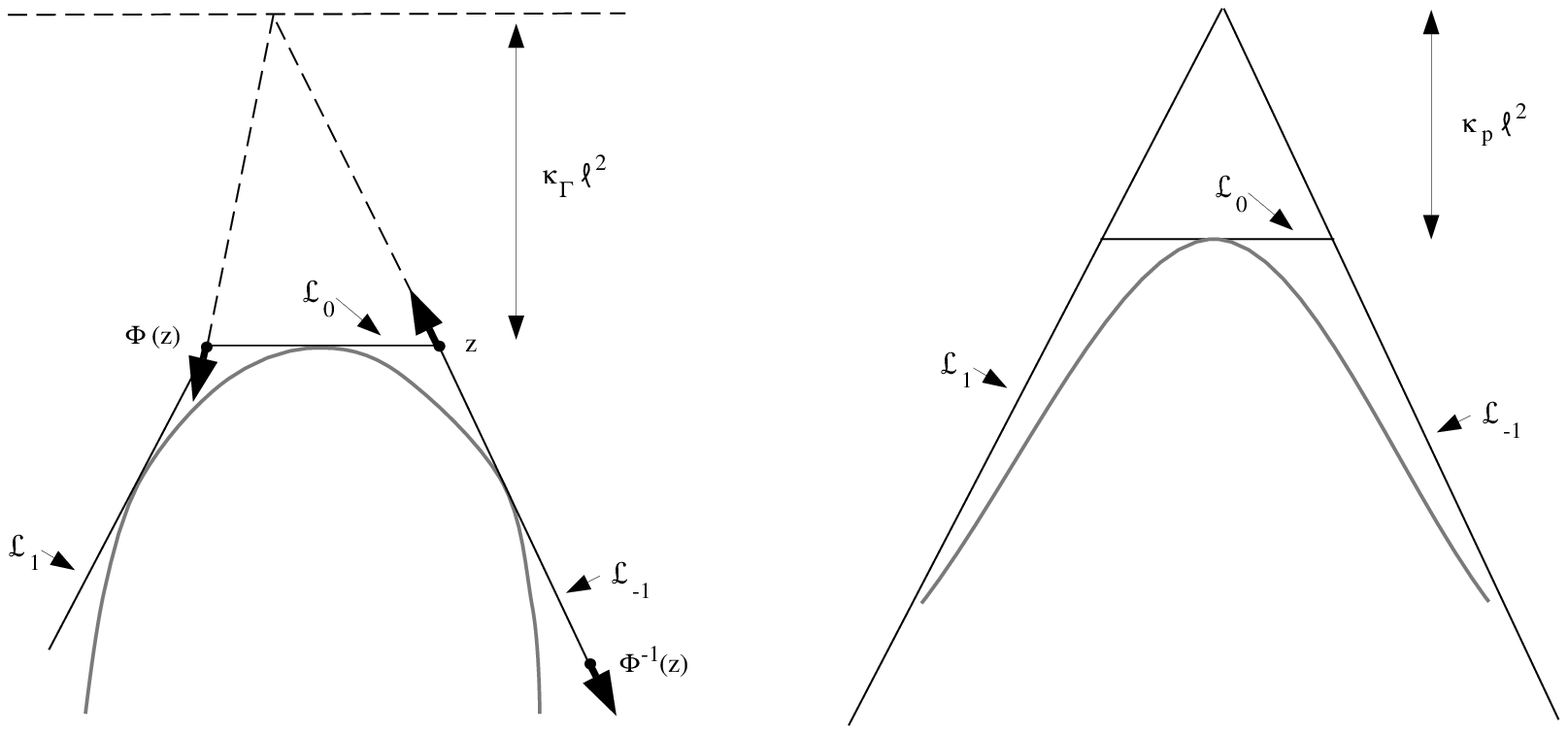,height=.45\hsize}{6.3}{(a) Tumbling tangents in 
dual billiards; (b) The comparison hyperbola}

It is geometrically clear now that when the 
curvature is infinite at the point of tangency of
 $\LL_0$ (\ie\ the radius of curvature 
vanishes) one always gets tumbling tangents for points $z$
near \mG, and thus no invariant circles near \mG.
It is also clear why the same argument does not work for
$z$ that are far away from \mG.

To get a more refined criterion, construct a hyperbola
as the area envelope inside the pair of lines
$\LL\mo$ and $\LL_1$ using an area equal to that of the
triangle bounded by  the $\LL\mo, \LL_0$, and
$\LL_1$ (see Figure 6.3(b)). This hyperbola will be tangent to $\LL_0$ at the
same point that \mG\ is.
 Now if $\Psi$ is the
\DB\ map of this hyperbola, the image of tangent vector based
at $z$ that points in the direction of $\LL\mo$
will be a vector based at $\Phi(z)$ that points in the direction
of $\LL_1$. Thus by \S 6.1, the point
$\LL\mo \cap \LL_1$ is a distance  $\kappa_p \ell^2$ from
$\LL_0$ where $\kappa_p$ is the curvature of the hyperbola 
at the point of tangency. Recalling the geometric
condition for a tumbling tangent for $\Phi$ one sees that
the tangent tumbles if and only if 
the bounce line for $\Gamma$ is higher than the bounce line
for the hyperbola, \ie\ if  $\kappa_\G > \kappa_p$. A simple
calculation shows that this condition is identical to (4.1).

\hfill\eject

\noindent{\bf \S 7 Dual billiards and impact oscillators}

In this section we show that the dynamics of a \DB\ map associated
with a convex curve is identical to the dynamics of a certain impact
oscillator. As in the previous section, our 
purpose is primarily descriptive so we
maintain a somewhat informal tone.

\medskip
\noindent{\bf \S 7.1 The impact oscillator.}
The first ingredient in the impact oscillator is a periodically
moving wall.
Its position at time $t$ is given by $p(t)$, and we
require that $p(t + 2 \pi) = p(t)$, $p$ is continuous, and
$p$ is sub-sine as defined in Remark 2.1.1.
 As noted there, if $p$ is twice differentiable this happens
if and only if 
$${\ddot p} + p := \rho \ge 0.\eqno(7.1)
$$ Alternatively, we 
can start with a periodic
 $\rho$ with $\int \rho(t) \exp(i t)\;dt
= 0$, pick initial conditions for $p(0)$ and $\dot p(0)$,
and then solve (7.1) for $p$.

The second ingredient in the impact oscillator is a particle
whose position is given by $x(t)$. Except at collision
this particle is a simple
harmonic oscillator with period $2 \pi$, \ie\ it satisfies
$$
{\ddot x} + x = 0.\eqno(7.2)
$$

The collisions of the particle with the wall are assumed to be
perfectly elastic and the wall has infinite mass; so
if $t_c$ is a time when $x(t_c) = p(t_c)$, then 
$$
-\bigl(\dot x_{before}(t_{c} ) -\dot p(t_{c})\bigr)= 
\dot x_{after}(t_{c}) - \dot p(t_{c}). \eqno(7.3)
$$

A {\de solution} of the impact oscillator is a function
$x(t) \ge p(t)$ that satisfies (7.2) except when 
$x(t) = p(t)$ (\ie\ collisions) at which time one starts a new solution to
(7.2) using the initial conditions dictated by (7.3).
We also require that $\dot p > \dot x$
at collision.  Thus if $\rho$ vanishes on a nontrivial
interval $J$, we do not consider an $x(t)$ with
$x(t) = p(t)$ for $t\in J$ a solution. Note that
since $p$ is sub-sine, if $\dot p > \dot x$
at one collision the same inequality holds at the next
collision, \ie\ the situation in Figure 7.1
{\it does not} occur. From a more physical
point of view, the fact that $p$ is sub-sine
implies that there are no so-called graze solutions
([BD]). In  these solutions the particle  
approaches the wall, touches
it for an instance matching the wall's position and velocity, and then pulls
away without a collision. The absence of these solutions
implies that the return time to the wall is a continuous function 
of $\dot x(t_c)$ and $t_c$.  

%\fig{fig7_1.ps,height=.30\hsize}{7.1}{A solution to the impact oscillator}

\fig{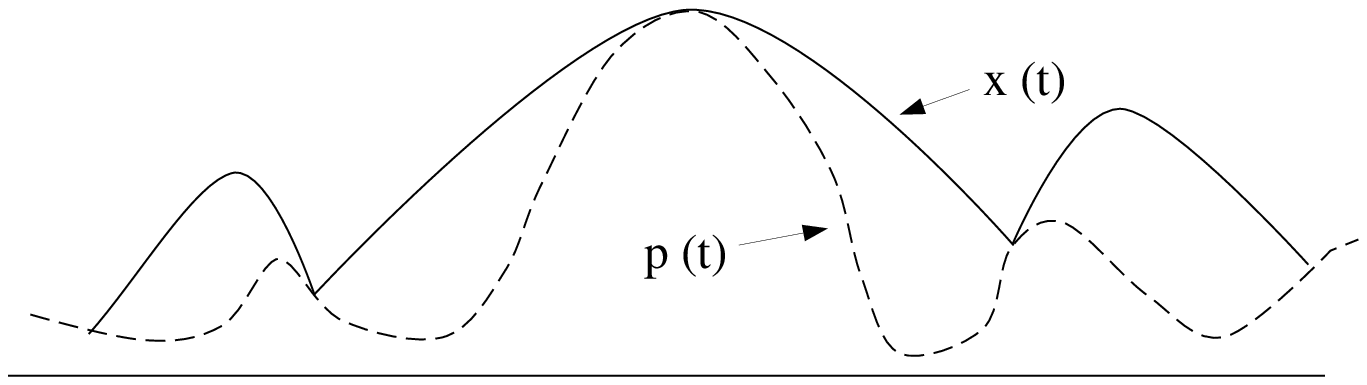,height=.20\hsize}{7.1}{Solution ruled out by sub-sine hypothesis}

\medskip
\noindent{\bf\S 7.2 Pictorial connection of the two systems.}
Fix a convex curve \mG\ in the plane with radius of curvature
function $\rho$, and assume for simplicity that  the origin is interior to \mG.
Call the orthogonal projection of \mG\ onto the $x$-axis  $I_0$ 
(see Figure 7.2(a)).
The distance of the origin to the right endpoint of $I_0$ is
$p(0)$ as defined in \S 2.1. Now imagine rotating the 
plane clockwise about the origin while keeping the direction of the 
$x$-axis fixed. If $I_\t$ is the projection of \mG\ onto the x-axis after
we have rotated by an angle $\t$, then the distance from the origin
to the right endpoint of $I_\t$ is $p(\t)$. We will
think of the segments $I_\t$ as the shadow of \mG\ and we will rotate
at unit speed and identify the angle $\t$ with the time $t$.
It follows from (2.3) that the
right endpoint of the shadow is evolving according
to $\ddot p + p = \rho$. Thus the end of
the shadow is moving like the wall in the impact oscillator.

\longfig{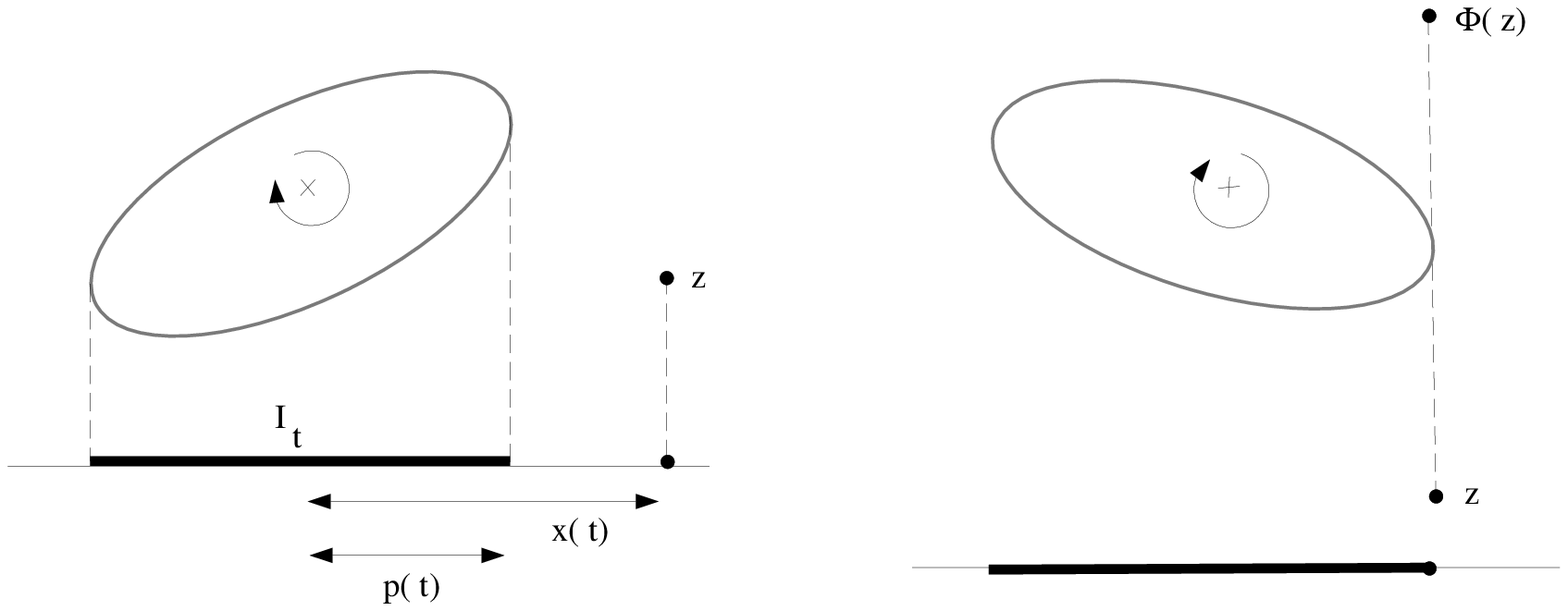,height=.40\hsize}{7.2}{(a)Pre-collision in 
the rotating plane \hfill\break model;
 (b)Collisions in the rotating plane model}

To introduce the particle, pick a point $z_0$ exterior
to \mG\ and monitor the evolution of its shadow (its orthogonal projection
onto the $x$-axis) as we rotate the plane. The point  moves in a 
circle and so if its shadow is given by $x(t)$, then
$\ddot x + x = 0$. 

Now for the collisions. Assume that  at time
$t$,  $x(t) > p(t)$ as in Figure 7.2(a). As we continue rotating,
at some time $t_c$, $x(t_c)$ will equal $p(t_c)$, \ie\ the shadow
of the point and that of \mG\ will collide. 
This will be when the supporting line to \mG\ that contains $z_0$ is vertical.
Since the shadows of  $z_0$ and \mG\ are 
taking the role of the particle and the wall
in the impact oscillator,  when the shadows collide
we should act on the point $z_0$ in such a way
that its shadow rebounds
from the shadow of \mG\ with the appropriate velocity.
It is clear from Figure 7.2(b), that the correct action
is to reflect
$z_0$ around the point of tangency to the \mG\ to obtain
a new point $\Phi(z_0)$. As the plane continues to rotate, the shadow
of $\Phi(z_0)$ will evolve away from the shadow of
\mG\ as if it had a perfectly elastic collision.
The map $\Phi$ just defined is clearly 
identical to the \DB\ map associated with \mG.

\medskip
\noindent {\bf  \S 7.3 Connecting the generating functions.}
To get a more precise connection between \DB\ and the impact
oscillator we connect the generating function of the 
\HTM\ associated with \DB\ with an action integral of an impact 
oscillator equivalent to  that of \S 7.1.

First we derive a new formula for the generating function of \DB. 
If $B(\t_1, \t_2)$ is the area of the sector of
\mG\ between the $\t_1$ and $\t_2$ (see Figure 7.3), then using
Stokes theorem 
$$\eqalign{ B(\t_1, \t_2) &=
  {1\over 2} \int_{\t_1}^{\t_2} \alpha(\t)\wedge {\bf n}_\t \; d s(\t)\cr
&={1\over 2} \int_{\t_1}^{\t_2} p(\t) \rho(\t)\; d\t.\cr}
$$
Recalling that $\rho = p + p^{\prime\prime}$, integrating by parts yields
$$B(\t_1, \t_2) =
  {1\over 2}\bigl( \int_{\t_1}^{\t_2} p^2 - (p^\prime)^2\; d\t
+ p(\t_2) p^\prime(\t_2) - p(\t_1) p^\prime(\t_1)\bigr).\eqno(7.4)$$

\fig{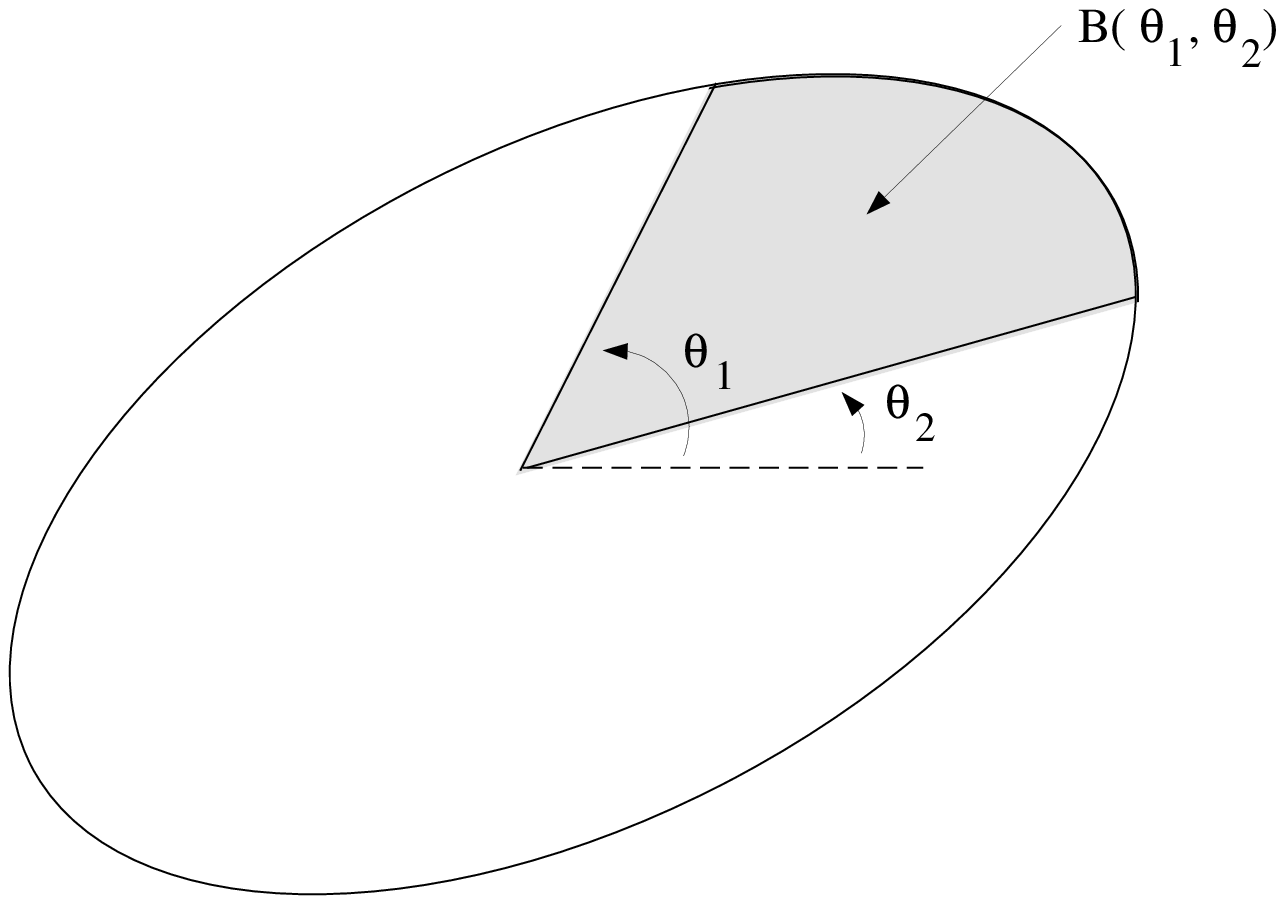,height=.30\hsize}{7.3}{The region $B(\t_1, \t_2)$}

We now move the origin to the point $Q= \LL_{\t_1}\cap \LL_{\t_2}$
and let $\eta(\t)$ be the distance from $Q$ to
$\LL_\t$.
With this new origin the area computed as 
 $B(\t_1, \t_2)$ is the 
generating function of the \HTM\ associated with the
\DB\ map (see Remark 3.2.1). In a slight abuse of notation 
we will treat this as a
function of $\t$ rather than $x$ and denote it $h_\G(\t_1, \t_2)$. 
To use (7.4) we note that from the point of
view of the point $Q$, the radius of curvature of
\mG\ is $-\rho(\t)$, and so using the 
fact that $\eta(\t_1) = \eta(\t_2) = 0$, we get
$$h_\G(\t_1, \t_2) =  
{1\over 2} \int_{\t_1}^{\t_2} 
(\eta^\prime)^2 - \eta^2  \; d\t.\eqno(7.5)
$$
where the function $\eta(\t)$ is the (unique!) solution to
$$\eta^{\prime\prime} + \eta = - \rho$$
with initial conditions $\eta(\t_1) = \eta(\t_2) = 0$.
For future reference note that from the geometric
interpretation of the derivative of the height function (see
\S 2.1), $\eta^\prime(\t_1) = R(\t_1, \t_2)$ and 
$\eta^\prime(\t_2) = -L(\t_1, \t_2)$. 

To connect this action to  the impact oscillator 
above we define for that system a new variable
 $\eta(t) = x(t) - p(t)$; thus 
$\eta$ measures the distance between the particle and the wall.
 The collision rule is now
$$-\deta_{before}(t_{c}) 
= \deta_{after}(t_{c}),\eqno(7.6)
$$ 
and motion away from collision is governed by 
$$\ddot \eta + \eta = -\rho.\eqno(7.7)$$
It will be useful to think of (7.6) and (7.7) as defining
a new impact oscillator that is dynamically equivalent
to the original one. 
In this new system, $\eta$ is the position of
a particle governed by  (7.7) that has perfectly  elastic collisions 
with a stationary wall.

Now between impacts this  system has Lagrangian 
$L(x, \dot x, t) := \dot x^2/2 -
x^2/2 - x\rho(t)$ and Hamiltonian $H(q, p, t)
:= q^2/2 +p^2/2 + p\rho(t)$. If $t_1 < t_2 < t_1 + \pi$, 
we define 
$$h_I(t_1, t_2) = \min \int_{t_1}^{t_2} L(x, \dot x, t)\;dt$$
where the minimum is taken over all paths $x(t)$ with
$x(t_1) = 0$ and $x(t_2) = 0$. Thus
$$h_I(t_1, t_2) = \int_{t_1}^{t_2} L(\eta, \deta, t)$$
where $\eta = \eta(t; t_1, t_2)$ is the (unique!) solution to 
$\ddot \eta + \eta = -\rho$ 
with boundary values $\eta(t_1) = 0 = \eta(t_2)$. 

Now note that $h_I$ is the Hamilton principal
function (see Remark 7.3.1) with the starting and 
ending configuration variables fixed at
zero, thus 
$$\eqalign{
\partial_1 h_I(t_1, t_2) &=H(\eta(t_1), \deta(t_1), t_1)  = \deta(t_1)^2\cr
\partial_2 h_I(t_1, t_2) &= -H(\eta(t_2), \deta(t_2), t_2)  = - \deta(t_2)^2.\cr}\eqno(7.8)$$
Recalling the collision rule
(7.6), we see that a sequence of times $(t_i)$ will be the sequence of
a collision times of a solution to the oscillator
precisely when it is gives stationary configuration 
as in \S 1.2
for the variational problem defined using
$h_I$.
But exactly as in (7.4) one has 
 $$\int_{t_1}^{t_2} \eta \rho = \int_{t_1}^{t_2}
\eta^2 - \deta^2,$$ and
so $h_\G(t_1, t_2) = - h_I(t_1, t_2)$. 

In conclusion, given an impact oscillator 
described by (7.6) and (7.7),  define 
$\Phi_I$ as  the return map to collisions, \ie\ 
$\Phi_I: (t_0, (\deta(t_0))^2) \mapsto (t_1, (\eta^\prime(t_1))^2)$ where
$t_0$ and $t_1$ are the times at consecutive collisions;
we have shown that $\Phi_I$ is identical to the \HTM\ map associated
with \DB\ on a convex curve with radius of curvature function 
$\rho$.

{\bf Remarks:}

{\bf 7.3.1:} If $H$ is a Hamiltonian that satisfies the
Legendre condition, then Hamilton's principal function is 
$$S(q_0, q_1, t_0, t_1) = \int_{t_0}^{t_1} L(x, \dot x,t)$$
where $L$ is the corresponding Lagrangian and 
$x(t)$ is a solution to Hamilton's equations (or equivalently an extremal
of the integral) 
with $x(t_0) = q_0$ and $x(t_1) = q_1$ (see \eg\ chapter 9 in 
[Gs]). For $S$ to be 
well-defined this solution needs to be unique.
The differential of the principal function
is $dS = p\;dq - H\; dt$. What was defined as 
$h_I$ above is the restriction of  $S$
to $0 = q_0 = q_1$. In [M5] Moser uses a similar construction,
but restricts $S$ to $t_0 = 0$ and $t_1 = 1$ to obtain the generating function
of the time-one map of the solution flow to a time-periodic Hamiltonian.
He shows that any twist map of the compact annulus can
be thus obtained. This raises the question
of which twist maps of the half-cylinder can arise 
from impact oscillators.

{\bf 7.3.2:}
One can also find solutions to the impact oscillator
by  considering  broken paths 
consisting of a collection of functions $x_i(t)$ and times
$t_i$ so that $x_i(t_i) = 0 = x_i(t_{i+1})$  and $x_i > 0$,
otherwise. Solutions are then extremals for
 $\sum L(x_i, \dot x_i, t)$. 

 In this context introducing the generating function $h_I$ is analogous to
the strategy of broken geodesics in the variational theory
of geodesics (see page 330 of [Gl] for some history).
 For extremals of the fixed endpoint problem the introduction of
$h_I$ replaces an infinite dimensional problem with a finite dimensional one.
Crucial to the application of this method is the
fact that that there is a unique extremal 
from $x=0$ to $x=0$ which begins at a time
$t_0$ and ends at $t_1$, where $t_0 < t_1 < t_0 + \pi$.

{\bf 7.3.3:} Yet another way of connecting the \DB\ map to
an impact oscillator can be obtained through a kind of
regularization of the collisions with the
wall. Solutions to  the equation
$$\ddot \eta + \eta = -\sgn(x) \rho(t)\eqno(7.9)$$
behave as if they pass through the wall instead of
colliding. Upon passage through the wall the sign of
the forcing is changed. 
Solutions of (7.9) for $x < 0$ are the reflections through the origin of those
for $x>0$. We can use (7.9) to generate a flow on
$\plane \times \circle$ via
$$\eqalign{
\dot x &= y\cr
\dot y & = -x - \sgn(x) \rho(\tau)\cr
\dot \tau &= 1\cr}$$
The return map of the flow to the open annulur cross section
defined by $x=0$, $ y > 0$ will be the second iterate
\DB\ map on the \mG\ with radius of curvature $\rho$.
But note that the map will be in $(\theta,\ell)$ coordinates
and will thus not be area-preserving.
 
The correct choice of second coordinate (\ie\ the one that
makes the map area-preserving) can be seen to be dictated by
(7.8) and Hamilton-Jacobi theory as in Remark 7.3.1.
From yet another point of view we
can make the $1 {1\over 2}$-degree
of freedom system given by (7.9) into an autonomous $2$-degree of freedom
solution in the usual way by letting
$H^\ast(q, p, E, \tau) = H(q, p, \tau) -  E$, where
now $H(q, p, t) = p^2/2 + q^2/2 + |q| \rho(t)$. 
Letting the $\tau$ variable be on a circle, and restricting attention to
$H^\ast = 0$, the return map in 
$(\tau, E)$ coordinates to $q=0, p>0$  will be the second iterate
of dual billiards map. Note that this map is area preserving because the
flow preserves volume and since $q=0$ on the chosen cross section,
the $E$ coordinate there is just $p^2/2$.

{\bf 7.3.4:} In [Kl], [SV] and [GS], it is shown
that \DB\ on rational polygons are stable, \ie\ there
do not exist any unbounded orbits.
Recalling \S 2.1, this implies that a class of 
impact oscillators as defined by (7.6) and (7.7)
whose forcing functions are contained in 
a certain collection of linear combinations
of delta functions have no unbounded solutions.

{\bf 7.3.5:} The physical content of the fact that the return map
to the wall  using the $(t, \deta^2)$ coordinates is a twist
map is simply that when  two points are started at the 
same time at the wall, the one with greatest initial velocity will take 
the longest to return to the wall.

\medskip
\noindent{\bf Appendix}

The appendix is a remark on invariant circles for (inner) billiards.
As is the case with dual billiards (\S 4.2), 
the solution to the inverse problem for the
existence of invariant circles is a special case of Proposition 1.3.
The remark rests on a elementary computation, 
but it does not seem to be widely known. Billiards is the subject
of a large body of work (\cf\ [T3]); we only briefly describe the system here.
Billiards as a twist map is considered in [M1], [D] and [Me2].

If \mG\ is an oriented  convex curve in the plane, the billiards map
is defined using a point particle moving freely in the region
bounded by \mG. The particle has perfectly elastic collisions
with the curve, so the angle of incoming motion with the tangent to the
curve is equal to the angle of the outgoing direction.
The billiards map just keeps track of  collisions with the curve.
The map can be described using coordinates $t$, the arc length parameter 
along the curve,  and $\omega$,
the angle between an outgoing line of motion and the tangent to the curve.
The billiards map is then $(t_0,\omega_0)\mapsto (t_1, \omega_1)$,
where  $(t_0,\omega_0)$ represents an initial point on the curve and an
outgoing direction, and
$(t_1, \omega_1)$ represents  the point on the curve at the 
next collision and its outgoing direction after the collision (see
Figure A.1).
This defines a homeomorphism of the annulus $\circle\times [0, \pi]$.

\fig{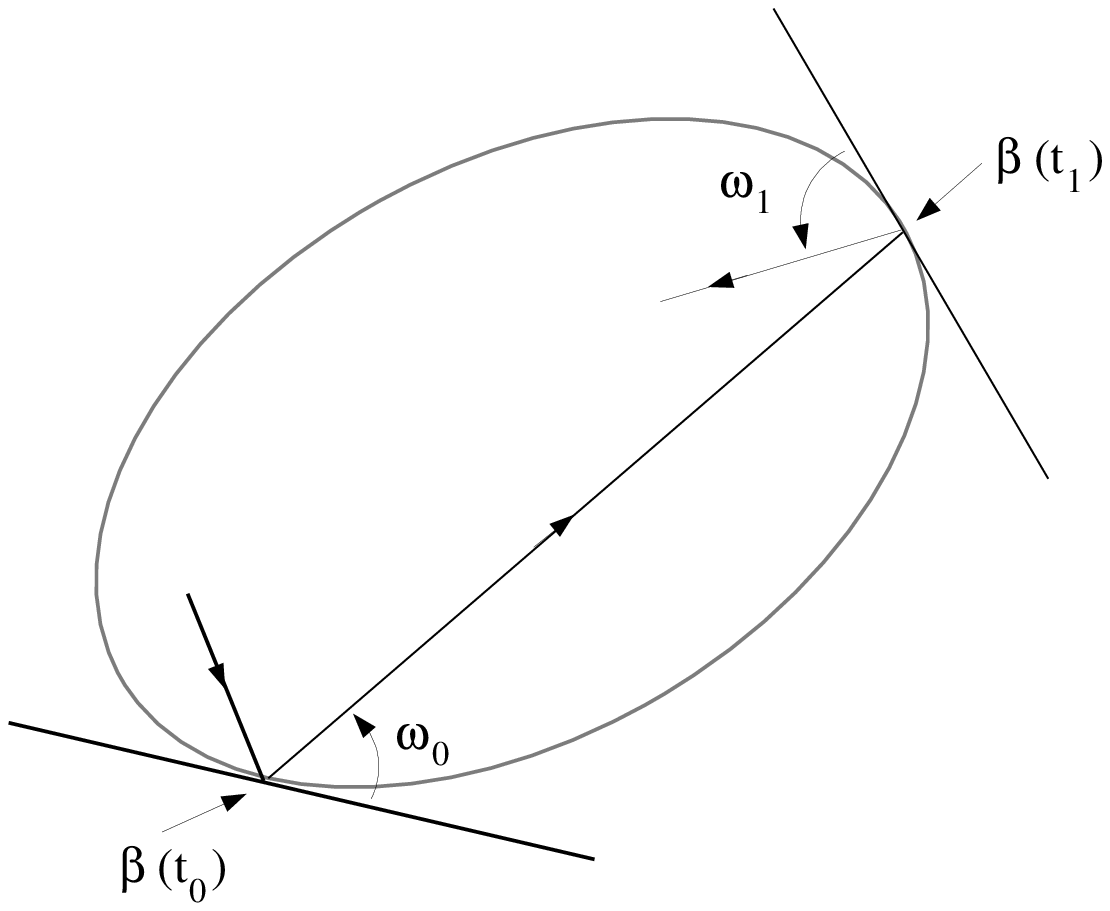,height=.35\hsize}{A.1}{The billiards system}

An invariant circle of this annulus map corresponds to a geometric object
called a caustic. For simplicity we restrict to the case of
convex caustics. A convex caustic is a convex curve inside of \mG\ 
with the property that an outgoing path that is tangent to $\Gamma_1$ on
one pass must also be tangent on passes after all future (and past) collisions.
The inverse problem for the existence of caustics is the following:
given a convex curve $\Gamma_1$, does there exist
another convex curve $\G$ so that
$\Gamma_1$ is a caustic for billiards on \mG? The solution
to this problem involves the string body construction. 
Fix the curve $\Gamma_1$ and a length of string that is
longer than the  perimeter of $\Gamma_1$. Loop the string around
$\Gamma_1$ and put a pencil in the loop, pull the string
tight and move the pencil
around $\Gamma_1$ keeping the string tight. The pencil will
draw out a curve $\G$ that is characterized by the constancy of
the length $B + C + D$ as shown in Figure A.2. It is
clear that \mG\ will have
$\Gamma_1$ as a caustic. The converse is also true; if \mG\ has
a caustic $\Gamma_1$, then \mG\ may be obtained as
a string body of $\G_1$. Put another way, given any $\G_2$ inside
\mG\ we can define a quantity called the {\de string parameter}
$L(t) = B(t) + C(t) + D(t)$ where 
the quantities $B$, $C$, and $D$ may vary with the point
$\beta(t)$. The curve $\Gamma_2$ is a caustic if and only if
$L(t)$ is  constant.
This result 
is attributed to Minasian in [S]; it also appears
in [P] (see also [Tr]). 

\fig{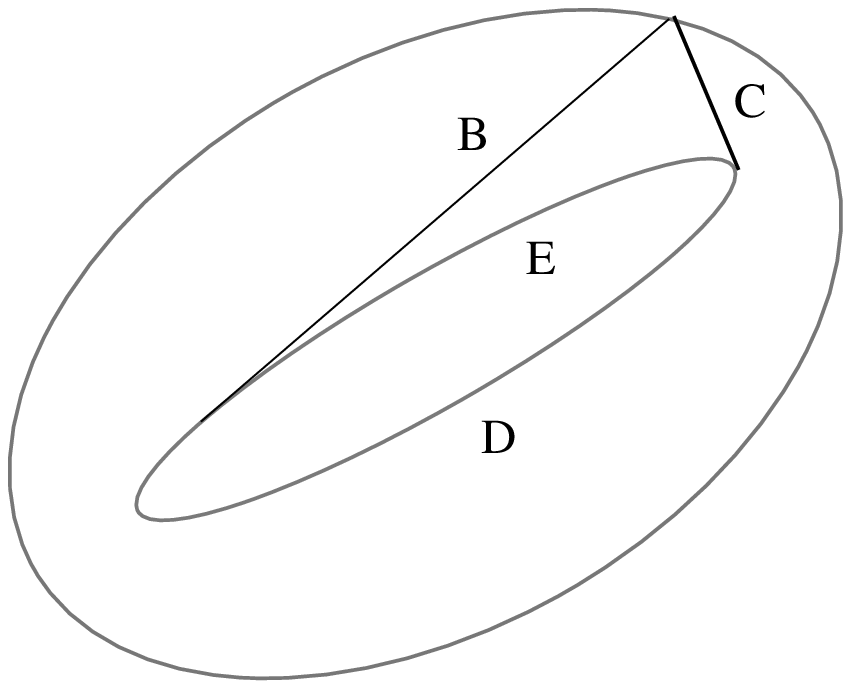,height=.30\hsize}{A.2}{The string parameter}

With the proper choice of coordinates, the annulus homeomorphisms
arising from billiards on  convex curves gives a class of
area preserving twist maps that have a theory similar to that
described in \S 1. In particular, there is a generating
function for the homeomorphism, and one can 
find orbits using  variational techniques. The analog of
Theorem 1.1 is true as well as much else. Our purpose here is
just to remark on one aspect of the theory, namely, the connection
of the area function $A(t)$ defined in \S 1.4 with the string body
construction. 

To  describe the billiards map more carefully,
begin by letting  $\beta(t)= (\beta_1(t), \beta_2(t))$ be a
parameterization of  convex curve \mG\ {\it by arc length}. 
For simplicity, assume that \mG\ has perimeter $2 \pi$, and so $t\in \circle$.
The generating function for the billiards map 
is the opposite of the length of the chord between impacts with the curve,  
$$h(t_0, t_1) = -\|\beta(t_1) - \beta(t_0)\|.$$
Thus, for example,
one can obtain  period three orbits by finding a 
triangle with edges on \mG\ that has
{\it maximal} perimeter among all such triangles.

In accord with (1.1) we differentiate $h$ to get the appropriate
second coordinates to make billiards an area-preserving twist map,
$$
h_1(t_0, t_1) = 
{\beta(t_1) - \beta(t_0) \over \|\beta(t_1) - \beta(t_0)\|} \cdot
\beta^\prime(t_0) = \cos(\omega(t_0, t_1))$$
where $\omega(t_0, t_1)$ is the angle labeled $\omega_0$ in  Figure A.1.
Thus the appropriate second coordinate is the opposite of
the cosine of  the angle between the outgoing ray and the tangent
line to the curve. Billiards defines a twist map
because $h_{12} < 0$. The billiards on
a convex curve \mG\  then defines an area-preserving
twist map  $\psi_\G:\circle\times [-1, 1]\ra\circle\times [-1, 1]$.	

Now assume that \mG\ has a caustic, \ie\ there exists
a $\psi_\G$-invariant circle.  
By the appropriate version of Birkhoff's Theorem (Theorem 1.1(b), see
[M1]), an invariant circle  
is the graph of a Lipschitz function $u:\circle\ra [-1,1]$.
The function whose graph is the lift of the invariant circle to
$\reals\times [-1,1]$ will also be denoted $u$.
The coordinates on both $\circle$ and $\reals$ will be called
$t$.

Examining the  geometric situation
one sees that a caustic $\G_1$ defines a one-parameter family of angles
each associated with a point on \mG. The angle
at the point $\beta(t)$ is denoted $\nu(t)$.  The point of tangency
to $\Gamma_1$ of the outgoing ray from $\beta(t)$ is denoted
$\sigma(t)$ and the advance map along \mG\ is
$g(t)$ (see Figure A.3).

\fig{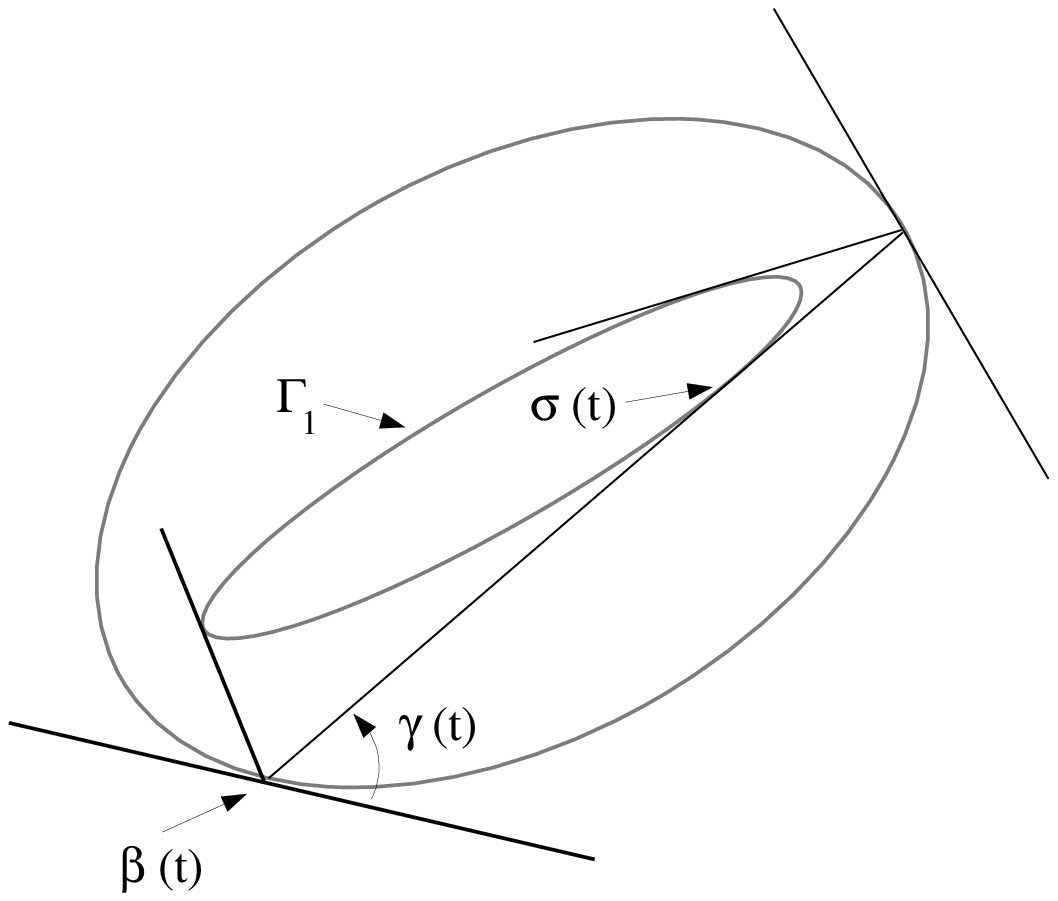,height=.35\hsize}{A.3}{Functions associated with a caustic}

Using the connection of $u(t)$ to these geometrically
defined quantities,
$$
u(t) = -\cos(\nu(t)) = - \beta^\prime(t) \cdot
{\sigma(t) - \beta(t) \over \|\sigma(t) - \beta(t)\|},
$$
and so 
$$
{d\over dt}(\|\sigma(t) - \beta(t)\|) = u(t) + 
{\sigma(t) - \beta(t) \over \|\sigma(t) - \beta(t)\|} \cdot\sigma^\prime(t).
 $$
The area function (\S 1.4) associated with the invariant circle   is 
$$\eqalign{A(t) &= \int_t^{g(t)} u(s) \; ds - h(t, g(t))\cr
&=-\int_t^{g(t)} {\sigma(s) - \beta(s) \over \|\sigma(s) - \beta(s)\|}
\cdot \sigma^\prime(s)\;ds + 
\int_t^{g(t)} {d\over ds}(\|\sigma(s) - \beta(s)\|)\; ds - h(t, g(t))\cr
&= -\int_t^{g(t)} \|\sigma^\prime(s)\|\; ds + \|\sigma(g(t)) - \beta(g(t))\| \cr
&\qquad\qquad - \|\sigma(t) - \beta(t)\| + \|\beta(g(t)) - \beta(t)\|\cr}$$
where we used the fact that 
${\sigma(s) - \beta(s) \over \|\sigma(s) - \beta(s)\|}$
is the unit vector in the direction of $\sigma^\prime(s)$.

The first term in the final expression for $A(t)$
is the opposite of the length along $\Gamma_1$ between the
two tangent rays (this length is labeled E in Figure A.2).
If the perimeter of $\Gamma_1$ is $P$,
this first  term is thus $D-P$. The second and third terms in
the expression for $A(t)$ are B and C as shown in Figure
A.2,  and so $A(t) = L(t) - P$.
Thus using the appropriate version of Proposition 1.3,
the string parameter is
constant if and only if $\G_1$ is a  caustic for \mG.

{\bf Remark A.1:} Note that in the billiards systems 
the generating function  is always
negative. The geometric interpretation of the generating function
as an area  in the annulus (Remark 1.1.1)  is still valid, but now
the  area is  below the $t$-axis. The area function in the case
of billiards is a positive
number giving the area below the invariant circle and
between a vertical arc and its image.

\medskip
\noindent{\bf References}

\item{[A]} Angenent, S., Monotone recurence relations, their
Birkhoff orbits and topological entropy, {\it  Ergod. Th. \& Dynam. Sys.},
{\bf 10},  1990, 15--41.

\item{[Bg]} Bangert, V., Mather sets for twist maps and geodesics on
tori, {\it Dyn. Rep.}, {\bf 1}, 1988, 1--56.

\item{[B+H]} Boyland, P. and Hall, G. R., Invariant 
circles and the order structure of periodic
orbits in monotone twist maps, {\it Topology}, {\bf 26}, 21-36, 1987.

\item{[BD]} Budd, C. and Dux, F., Intermittency in impact oscillators close to
resonance, {\it Nonlinearity}, {\bf 7}, 1994, 1191--1224.

\item{[Dv]} Devaney, R.,{ \it An Introduction to Chaotic Dynamical Systems,}
Addison-Wesley, 1989.

\item{[D]} Douady, R., Th\`ese de Troisi\`eme Cycle, Univ. Paris 7, 1982.

\item{[E]} Eggleston, H., {\it Convexity}, Cambridge University Press,
1958.

\item{[FT]} Fuchs, D. and \Tab, S., Segments of equal area, 
{\it Quantum}, {\bf 2}, 1992.

\item{[Gs]} Goldstein, H., {\it Classical Mechanics}, second
edition, Addison--Wesley, 1980.

\item{[Gl]} Gol\'e, C., Periodic orbits for Hamilitonian systems
in cotangent bundles, {\it Trans. AMS}, {\bf 343}, 1994, 327--347.

\item{[G1]} Green, J., Sets subtending a constant angle on a circle,
{\it Duke Math. Jour.}, {\bf 17}, 1950, 263--267.

\item{[G2]} Green, J., Length and area of a convex curve under affine
transformation, {\it Pac. J. Math.}, {\bf 3}, 1953, 393--402.

\item{[GK]} Gutkin, E.  and Katok, A., Caustics for innner and outer
billiards, preprint, 1994.

\item{[GS]} Gutkin, E. and Simanyi, N., Dual polygonal billiards and necklace
dynamics, {\it Comm. Math. Phys.}, {\bf 143}, 1991, 431--450.

\item{[Ha]} Halpern, B., Strange billiard tables, {\it TAMS}, {\bf
232}, 1977, 297--305.

\item{[H]} Herman, M., Sur les courbes invariantes par les diff\'eomorphismes de
l'anneau (volume 1), {\it Ast\'erisque}, {\bf 103--104}, 1983.

\item{[Hr]} Hubacher, A., Instability of the boundary in the billiard
ball problem, {\it Comm. Math. Phys.}, {\bf 108}, 1987, 483--488. 

\item{[K]} Katok, A.,  Some remarks on the Birkhoff  and Mather twist
theorems, {\it Erg. Th. \& Dynam. Syst.}, {\bf 2}, 1982. 

\item{[Kl]} Kolodziej, R., The antibilliard outside a polygon, {\it Bull. Pol.
Acad. Sci}, {\bf 37}, 1989, 163--168.

\item{[Mc]} MacKay, R., Introduction to the dynamics of area presrving
maps, {\it Physics of Particle Accelerators}, ed. by Month, M. and Dienes,
M., Am. Inst. Phys. Conf. Proc., {\bf 153}, vol. 1, 1987, 534--602.

\item{[MP]} MacKay, R. and Percival, I.,  Converse KAM: Theory and practise,
{\it Comm. Math. Phys.}, {\bf 98}, 1985, 469--512.

\item{[M1]} Mather, J., Glancing billiards, {\it Ergod. Th. \& Dynam.
Sys.}, {\bf 2}, 1982, 397--403.

\item{[M2]} Mather, J., Nonexistence of invariant circles, 
{\it Ergod. Th. \& Dynam. Sys.},  {\bf 4}, 1982, 301--309.

\item{[M3]} Mather, J., Letter to R. MacKay, 1984.

\item{[M4]} Mather, J., Variational construction of orbits of
twist homeomorphisms, {\it Jour. A.M.S.}, {\bf 4}, 1991, 207--263.

\item{[M5]} Mather, J., Variational construction of connecting orbits,
{\it Ann. Inst. Fourier, Grenoble}, {43}, 1993, 1349--1386.

\item{[MH]} Mayer, K. and Hall, G. R., {\it Introduction to Hamilitonian 
Dynamical Systems and the N-body Problem}, Applied Mathematical
Sciences, {\bf 94}, Springer-Verlag, 1991.

\item{[Me1]} Meiss, J., Symplectic maps, variational principles, and transport,
{\it Rev.  Mod. Phys.}, {\bf 64}, 1992, 795--848.

\item{[Me2]} Meiss, J., Cantori for the stadium billiard, {\it Chaos},
{\bf 2}, 1992, 267--272.

\item{[Ms1]} Moser, J., Is the solar system stable?, {\it Math. Intell.}, {\bf
1}, 1978, 65--71.

\item{[Ms2]} Moser, J., {\it Stable and Random Motions in Dynamical Systems},
Ann. of Math. Stud., Princeton University Press, {\bf 77}, 1973.

\item{[Ms3]} Moser, J., Recent developments in the theory of 
Hamiltonian dynamical systems, {\it SIAM Rev.}, {\bf 28}, 1986, 459--485.

\item{[Ms4]} Moser, J., Breakdown of stability, {\it Nonlinear
Dynamics Aspaects of Particle Accelorators}, ed. by Jowett, J.,
Springer LNIP, {\bf 247}, 1986, 492--518.

\item{[Ms5]} Moser, J., Monotone twist maps and the calculus of
variations,
{\it Ergod. Th. \& Dynam. Sys.}, {\bf 6}, 1986, 401--413.

\item{[P]} Poritsky, H., The billiard ball problem on a table with convex
boundary -- a illustrative problem, {\it Ann. Math.}, {\bf 51}, 1950,
446--470.

\item{[SV]} Shaidenko, A. and Vivaldi, F., Global stability of a class
of discontinuous dual billiards, {\it Comm. Math. Phys.}, {\bf 110},
1987, 625--640.

\item{[S]} Sinai, Y., {\it Introduction to ergodic theory}, Princeton Univ.
Press, 1976. 

\item{[T1]} \Tab, S., Dual biliards, {\it Russ. Math. Surv.}, {\bf 48}, 1993,
75--102.

\item{[T2]} \Tab, S., On the dual billiard problem, {\it Adv. in Math.}, to
appear.

\item{[T3]} \Tab, S., Billiards, University of Arkansas Preprint  UofA-R-69, 1994.

\item{[Tr]} Turner, P., Convex caustics for billiards in $\reals^2$ 
and $\reals^3$, {\it Convexity and related combinatorial geometry}, 
Marcel Dekker, 1982, 85--105.

\item{[W]} Wojtkowski, M., Principles for the design of billiards with
nonvanishing Lyapunov exponents, {\it Comm. Math. Phys.}, {\bf 105},
1986, 391--441.

\bye